\documentclass[english,11pt]{article}

\usepackage[german,french,english]{babel}
\usepackage[cp850]{inputenc}
\usepackage{latexsym,graphicx, fancybox}
\usepackage{amssymb, amsmath, amsfonts}
\usepackage{color}
\usepackage{latexsym}

\font\tenmath=msbm10 scaled 1200
\font\sevenmath=msbm7 scaled 1200
\font\fivemath=msbm5 scaled 1200

\newcommand{\vertiii}[1]{{\left\vert\kern-0.25ex\left\vert\kern-0.25ex\left\vert #1 
    \right\vert\kern-0.25ex\right\vert\kern-0.25ex\right\vert}}

\newfam\mathfam \textfont\mathfam=\tenmath
\scriptfont\mathfam=\sevenmath \scriptscriptfont\mathfam=\fivemath
\def\math{\fam\mathfam}
\def\R{{\math R}}
\def\N{{\math N}}
\def\E{{\math E}}

\def\P{{\math P}}

\newtheorem{Remark}{Remark}[section]

\newtheorem{Thm}{Theorem}[section]
\newtheorem{Lem}{Lemma}[section]
\newtheorem{Pro}{Proposition}[section]

\newtheorem{Dfn}{Definition}[section]

\newfam\mathfam \textfont\mathfam=\tenmath
\scriptfont\mathfam=\sevenmath \scriptscriptfont\mathfam=\fivemath
\def\math{\fam\mathfam}

\def \^#1{\if#1i{\accent"5E\i}\else{\accent"5E#1}\fi}

\def \ind {1 \mkern -5mu \hbox{I}}

\def \cqfd{\quad_\Box}

\begin{document}
\selectlanguage{english}
\title{\bf Convex order, quantization and monotone approximations of {\em ARCH} models}
 
\author{ 
{\sc Benjamin Jourdain} \thanks{Universit\'e Paris-Est, Cermics (ENPC), INRIA, F-77455 Marne-la-Vall\'ee, France. E-mail: {\tt   benjamin.jourdain@enpc.fr}}
\and   
{\sc  Gilles Pag\`es} \thanks{Laboratoire de Probabilit\'es, Statistique et Mod\'elisation, UMR~8001, Campus Pierre et Marie Curie, Sorbonne Universit\'e case 158, 4, pl. Jussieu, F-75252 Paris Cedex 5, France. E-mail:{\tt  gilles.pages@upmc.fr}}}
\date{}
\maketitle 
\begin{abstract}We are interested in proposing approximations of a sequence of probability measures in the convex order by finitely supported probability measures still in the convex order. 
  We propose to alternate transitions according to a martingale Markov kernel mapping a probability measure in the sequence to the next and dual quantization steps. In the case of {\em ARCH} models and in particular of the Euler scheme of a driftless Brownian diffusion, the noise has to be truncated to enable the dual quantization step. We analyze the error between the original {\em ARCH} model and its approximation with truncated noise and exhibit conditions under which the latter is dominated by the former in the convex order at the level of sample-paths. Last, we analyze the error of the scheme combining the dual quantization steps with truncation of the noise according to primal quantization.
  \end{abstract}

\section{Introduction}\label{sec:1}
For $d\!\in\N^*$, and $\mu,\nu$ in the set ${\mathcal P}(\R^d)$ of probability measures on $\R^d$, we say that $\mu$ is smaller than $\nu$ in the convex order and denote $\mu\le_{cvx}\nu$ if
\begin{equation}
   \forall \varphi:\R^d\to\R\mbox{ convex },\;\int_{\R^d}\varphi(x)\mu(dx)\le\int_{\R^d}\varphi(y)\nu(dy),\label{eq:defcvx}
\end{equation}
when the integrals make sense (since any real valued convex function is bounded from below by an affine function $\int_{\R^d}\varphi(x)\mu(dx)$ makes sense in $\R\cup\{+\infty\}$ as soon as $\int_{\R^d}|x|\mu(dx)<+\infty$). We then also write $X\le_{cvx}Y$ for $X$ and $Y$ random vectors respectively distributed according to $\mu$ and $\nu$.

For $p\ge 1$, we denote by ${\mathcal P}_p(\R^d)=\{\mu\!\in{\mathcal P}(\R^d):\int_{\R^d}|x|^p\mu(dx)<+\infty\}$ the Wasserstein space with index $p$ over $\R^d$. When $\mu,\nu\in{\mathcal P}_1(\R^d)$, according to the Strassen theorem~\cite{Strassen},  $\mu\le_{cvx}\nu$ if and only if there exists a martingale coupling between $\mu$ and $\nu$ that is a probability measure $M(dx,dy)$ on $\R^d\times\R^d$ with marginals $\int_{y\in\R^d}M(dx,dy)$ and $\int_{x\in\R^d}M(dx,dy)$  equal to $\mu(dx)$ and $\nu(dy)$ respectively such that $M(dx,dy)=\mu(dx)m(x,dy)$ for some Markov kernel $m$ sharing the martingale property: $\forall x\!\in\R^d$, $\int_{\R^d}|y|\,m(x,dy)<+\infty$ and $\int_{\R^d}y\,m(x,dy)=x$. If $(X,Y)$ is distributed according to $M$, then $X$ and $Y$ are respectively distributed according to $\mu$ and $\nu$ and $\E\,(Y|X)=X$.

In this paper, we are interested in constructing approximations of a sequence $(\mu_k)_{k=0:n}\!\in{\mathcal P}_1(\R^d)^{1+n}$ in increasing convex order ($\forall\,k=0:n-1,\;\mu_k\le_{cvx}\mu_{k+1}$) by a sequence $(\widehat\mu_k)_{k=0:n}$ of probability measures with finite supports still in the convex order. 

An important   motivation comes from mathematical finance in order to provide model free bounds for exotic options written on discrete time path-dependent payoffs when the quoted prices of vanilla options written on the underlying asset(s) expiring at these discrete times are available and an absence of arbitrage opportunities assumption is made on the underlying asset(s) dynamics. This approach was formalized 
by the Martingale Optimal Transport (MOT) problem introduced in~\cite{Beietal} and  has received  recently a great attention in the financial mathematics literature. Thus, the structure of martingale optimal transport couplings~\cite{BeJu,CampLaachMart,DeMarchTouzi,Ghousskimlim,HLTo}, continuous time formulations~\cite{DoSo,GalHTo,Toetal}, links with the Skorokhod embedding problem~\cite{Beicoxhues}, numerical methods~\cite{AlJo,AlJo2,DeMarch,GuOb,HL19} and stability properties~\cite{BackPam,JouMar,Wiesel} have been investigated. The dynamics of underlying assets is often modeled by (martingale) diffusion processes. Such processes are  classically time discretized by their  Euler scheme, which belong to the wide class {\em ARCH} models. On the other hand, the two instants MOT problem between finitely supported distributions can be reduced to a numerically solvable  constrained optimization problem as emphasized in~\cite{AlJo2}. These two facts  lead naturally to tackle the problem of space discretization of probability distributions that preserve convex order. The main objective of this paper is  to propose a general method to solve this two-fold problem: spatially  discretize martingale {\em ARCH models} in such a way that the resulting marginals are finitely supported and preserve their martingale property. By ``finitely'', we underline with controllable size from a computational viewpoint.

However, to our  best knowledge, few studies consider the problem of preserving the convex order while approximating a sequence of probability measures, like the marginals of a discrete time process. We mention the thesis of Baker~\cite{Baker} who proposes the following construction in dimension $d=1$. Let for $u \in (0,1)$, $F_{k}^{-1}(u)=\inf \{x \in \R: \mu_{k}((-\infty,x])\ge u \}$  be the quantile of $\mu_k$ of order $u$. Let $(N_k)_{k=0:n}$ be a sequence of elements of $\N^*$ such that for $k=0:n-1$, $N_{k+1}/N_k\in\N^*$.
One has $\widehat{\mu}_0\le_{cvx}\widehat{\mu}_1\le_{cvx}\hdots\le_{cvx}\widehat{\mu}_n$ for the choice 
$$
\widehat{\mu}_k=\frac{1}{N_k}\sum_{i=1}^{N_k}\delta_{N_k \int_{\frac{i-1}{N_k}}^{\frac{i}{N_k} } F_k^{-1}(u) du },\;k=0, \ldots,n.
$$
{\em Dual (or Delaunay) quantization} introduced by Pag\`es and Wilbertz~\cite{PaWi0} and further studied in~\cite{PaWi1,PaWi2, PaWi3} gives another way to preserve the convex order in dimension $d=1$ (see the remark after Proposition~10 in~\cite{PaWi1}) when $\mu_n$ is compactly supported.

In two recent papers~\cite{AlJo,AlJo2}, Alfonsi, Corbetta and Jourdain propose to restore the convex ordering from any finitely supported approximation $(\widetilde \mu_k)_{k=0:n}$ of $(\mu_k)_{k=0:n}$. In dimension $d=1$, one may define the increasing (resp. decreasing) convex order by adding the constraint that the test function $\varphi$ is non-decreasing (resp. non-increasing) in~\eqref{eq:defcvx}. Moreover, according to~\cite{AlJo2}, this can be performed by forward (resp. backward) induction on $k$ by setting $\widehat{\mu}_0=\widetilde\mu_0$ (resp. $\widehat{\mu}_n=\widetilde\mu_n$) and computing $\widehat{\mu}_k$ as the supremum between $\widehat\mu_{k-1}$ (resp. infimum between $\widehat\mu_{k+1}$) and $\widetilde\mu_k$ for the increasing convex order when $\int_\R x\widetilde\mu_k(dx)\le \int_\R x \widehat\mu_{k-1}(dx)$ (resp. $\int_\R x\widetilde\mu_k(dx)\ge \int_\R x \widehat\mu_{k+1}(dx)$) and the decreasing convex order when $\int_\R x\widetilde\mu_k(dx)\ge \int_\R x \widehat\mu_{k-1}(dx)$ (resp. $\int_\R x\widetilde\mu_k(dx)\le \int_\R x \widehat\mu_{k+1}(dx)$). For a general dimension $d$,~\cite{AlJo} suggests to set $\widehat\mu_n=\widetilde\mu_n$ and compute by backward induction on $k=0:n-1$, $\widehat\mu_k$ as the projection of $\widetilde \mu_k$ on the set of probability measures dominated by $\widehat\mu_{k+1}$ for the quadratic Wasserstein distance by solving a quadratic optimization problem with linear constraints.

For general dimensions $d$ but with only two marginals ($n=1$) and $\mu_1$ compactly supported, the convex order is preserved by defining $\widehat{\mu}_0$ as a stationary primal (or Voronoi) quantization of $\mu_0$ on $N_0$ points and $\widehat{\mu}_1$ as a dual (or Delaunay) quantization of $\mu_1$ on $N_1$ points. 
In a companion paper~\cite{JP20a} more focused on martingale couplings and the Martingale Optimal Weak Transport problem which is a generalization of MOT introduced in \cite{BackPam}, we prove that when these quantizations are optimal and $N_0$ and $N_1$ go to infinity, then the value function of the MOWT problem for the marginals $\widehat{\mu}_0$ and $\widehat{\mu}_1$ converges to the value function for the marginals $\mu_0$ and $\mu_1$ under mild regularity assumptions on the cost function.

By Strassen theorem, there exists a sequence $(m_k)_{k=0:n-1}$ of martingale Markov kernels on $\R^d$ such that for each $k=0:n-1$, $\int_{x\in\R^d}m_k(x,dy)\mu_k(dx)=\mu_{k+1}(dy)$. For such a discrete time family of Markov transitions $(m_k(x,dy))_{x \in \R^d}$ satisfying the above equations  there exist measurable functions $(G_k:\R^d\times\R^q\to\R^d)_{k=0:n-1}$ and independent $\R^q$-valued random vectors $(Z_k)_{k=1:n}$ independent from $X_0$ distributed according to $\mu_0$ (see e.g. Lemma 2.22 in~\cite{Kallenberg} for the particular case $q=1$ and $Z_k$ uniformly distributed on $[0,1]$ for each $k=1:n$), such for each $k=0:n-1$ and each $x\in\R^d$, $G_k(x,Z_{k+1})\sim m_k(x,.)$. Then the sequence $(X_k)_{k=0:n}$ defined inductively by 
\begin{equation}\label{eq:Markovchain}
X_{k+1}= G_k(X_k,Z_{k+1}), \quad k=0:n-1,
\end{equation}
is a martingale and a Markov chain with transition kernels $(m_k)_{k=0:n-1}$ and for each $k=0:n$, $X_k\sim\mu_k$.

Have in mind that, as a consequence of the general stationarity property of dual (or Delaunay) quantization,  for any random vector  $(Y,U)\!\in L_{\R^{d+1}}^{\infty}(\Omega,{\cal A}, \P)$, $Y\perp\!\!\!\perp U$, $U\stackrel{d}{=}U([0,1])$, and any grid $\Gamma\subset\R^d$, $\widehat Y^{\Gamma, del}= {\rm Proj}_{\Gamma}^{del}(Y,U)$ satisfies
\[
\E\big(\widehat Y^{\Gamma, del}\,|\,Y \big) = Y,
\]
where  ${\rm Proj}_{\Gamma}^{del}$  is defined in Appendix~\ref{App:A2} (see~\eqref{eq:ProjDel}).

In order to dually quantize the chain $(X_k)_{k=0:n}$ we need exogenous i.i.d. random variables $U_k\sim U([0,1])$, $k=1:n$,  independent of $(Z_k)_{k=1:n}$ and $X_0$. In this {\em dual quantization} of the chain note that the starting value $X_0$ will have a special status and does not need to be dually quantized but simply spatially discretized by any tractable mean. A natural choice is then to perform on $X_0$ a Voronoi (primal)  quantization.  
%

\medskip
For every $k=0, \ldots,n$, we consider grids $\Gamma_k \subset \R^d$ satisfying the following inductive property:
\[
({\cal G}_{\Gamma,Z})\; \equiv\quad \forall\, x\!\in \Gamma_{k-1}, \quad G_{k-1}\big(x,Z_k(\Omega)\big)\subset {\rm conv}\big(\Gamma_k\big), \quad  k= 1:n.
\]
It is clear that, by induction on $k$, such $n+1$-tuples of grids satisfying $(\mathcal{G}_{\Gamma,Z})$ exist as soon as the mappings $G_k$ satisfy
\begin{equation}\label{eq:HypoBounded}
(\mathcal{G}_{G,Z})\; \equiv\quad \forall\,k=1:n,\;\forall\, x\!\in\R^d,\;  \quad G_{k-1}\big(x, Z_k(\Omega)\big) \mbox{ is bounded in }\R^d.
\end{equation}

 Then, we may define by induction
\begin{align}
\label{eq:def0} \widehat X_0 & = g_0(X_0)\mbox{ for some Borel function }g_0:\R^d\to\Gamma_0
                               ,\\
\label{eq:def2} \widetilde X_{k}&= G_{k-1}(\widehat X_{k-1}, Z_{k})\quad \mbox{ and }\quad\widehat X_k = {\rm Proj}_{\Gamma_k}^{del}(\widetilde X_k,U_k), \quad  k=1:n,
\end{align}
%

\noindent where  $(U_k)_{k=1:n}$ is  an i.i.d. uniformly distributed on the unit interval sequence,  independent of $(Z_k)_{k=1:n}$ and $X_0$.

The sequence $(\widehat X_0,\widetilde X_1,\widehat X_1,\cdots,\widetilde X_n,\widehat X_n)$ is then a martingale Markov chain with respect to the filtration $(\sigma(X_0),\sigma(X_0,Z_1),\sigma(X_0,Z_1,U_1),\cdots,\sigma(X_0,(Z_\ell,U_\ell)_{\ell=1:n-1},Z_n),\sigma(X_0,(Z_\ell,U_\ell)_{\ell=1:n}))$, and in particular $$\widehat X_0\le_{cvx} \widetilde X_1\le_{cvx}\widehat X_1\le_{cvx} \cdots\le_{cvx} \widetilde X_n\le_{cvx} \widehat X_n.$$
The chain $(\widehat X_k)_{k=0:n}$ is then an approximation of the chain $(X_k)_{k=0:n}$ (see Remark \ref{remconvgen} further on  for quadratic error estimations under Lipschitz continuity of the functions $G_k$ with respect to their first variable) with each $\widehat X_k$ taking finitely many values and the distributions $\widehat\mu_k$ of the random variables $\widehat X_k$ are finitely supported approximations of the $\mu_k$ such that $\widehat\mu_0\le_{cvx}\widehat\mu_1\le_{cvx}\cdots\le_{cvx}\widehat \mu_n$. Notice that in order to obtain a practical numerical procedure, one should also replace the noise $(Z_k)_{k=1:n}$ by an approximation taking finitely many values. But for general functions $G_k$, it is not clear at all how to construct an approximation preserving the martingale property.

The present paper is specialized to {\em ARCH} models 
\begin{equation}\label{eq:ARCH0}
X_{k+1}=X_k+\vartheta_k(X_k)Z_{k+1},\quad k=0:n-1.
\end{equation}
with $(Z_k)_{k=1:n}$ an $\R^q$-valued {\em white noise} (\footnote{By {\em white noise} we mean here  a sequence of independent square integrable  centered random vectors with identity covariance matrix.}) independent of $X_0$ and, for $k=0, \ldots,n-1$, $\vartheta_k$ goes from $\R^d$ to the space ${\mathbb M}_{d,q}$ of real matrices with $d$ rows and $q$ columns. We will consider possibly unbounded ``innovations" $(Z_k)_{k=1:n}$, with in mind applications to the Euler scheme of martingale Brownian diffusions.

The multiplicative structure of the function $G_k(x,z)=\vartheta_k(x)z$ has the following nice consequence : the martingale property is preserved when replacing the noise $(Z_k)_{k=1:n}$ by an approximation $(\breve Z_k)_{k=1:n}$ possibly taking finitely many values as soon as the random vectors $\breve Z_k$ are centered, mutually independent and independent of $X_0$. 
The remaining problem, especially when in presence of several times steps ($n\ge 2$), is to control at every time $k$ the (finite) size of the support of the approximation of $X_k$ : if  $X_0=x_0\!\in \R^d$ and   the $Z_k$ are  replaced by $\breve Z_k$ taking e.g. $N$ values, then $X_n$ will take $N^n$ values which is totally unrealistic as soon as\dots $N=2$ if $n=20$. The interest of performing a dual quantization step at each time is then to reduce the size of the support of the approximation of each $X_k$, $k=1:n$ while preserving the martingale property. 

Section \ref{sec:scheme} is devoted to the scheme obtained by replacing the noise by some bounded approximation possibly obtained trough primal quantization and then performing a dual quantization step at each time. We analyse the quadratic error of this scheme.

In section \ref{sec:convex}, we give conditions on the functions $(\vartheta_k)_{k=0:n-1}$ ensuring the convex ordering at the level of paths between the {\em ARCH} model with bounded approximate noise on the one hand and the original {\em ARCH} $(X_k)_{k=0:n}$ or the above scheme on the other hand.

Last, in Section \ref{sec:euler}, we make the previous results more specific in the important special case when the original {\em ARCH} model is the Euler scheme of a martingale Brownian diffusion.

\medskip
%
\noindent {\sc Definitions and notations.}

\noindent $\bullet$ The space of real matrices with $d$ rows and $q$ columns is denoted by ${\mathbb M}_{d,q}$. 

\noindent $\bullet$ $|\cdot|$ denotes the canonical Euclidean norm on $\R^d$.



\noindent $\bullet$ ${\rm conv}(A)$ denotes the (closed) convex hull of $A\subset \R^d$ and ${\rm card}(A)$ or $|A|$  its cardinality (depending on the context).

\noindent $\bullet$  For $p\!\in [1,+\infty)$, let $$\mathcal{W}_p(\mu, \nu)=\inf\Big\{\left(\int_{\R^d\times\R^d}|x-y|^p M(dx,dy)\right)^{1/p},\; \mbox{ where $M$ has marginals }\mu\mbox{ and }\nu\Big\}\le +\infty$$ denote the Wasserstein distance with index $p$. This is a complete metric on the set $\mathcal{P}_p(\R^d)$ of probability distributions on $\big(\R^d, {\cal B}or(\R^d)\big)$ with finite $p$th moment.

\noindent $\bullet$ For every integer $N\ge 1$, we denote by  $\mathcal{P}(\R^d, N)$ the set of distributions on $\R^d$ whose support contains at most $N$ points. 

\noindent $\bullet$ The symbol $\perp\!\!\!\perp$ denotes the independence of random variables or vectors.

\noindent $\bullet$ A white noise is a sequence of independent square integrable  centered $\R^q$-valued random vectors with identity covariance matrix $I_q$.

\section{Quantized  martingale approximation scheme of an {\em ARCH} model}\label{sec:scheme}
\subsection{The quantized scheme}\label{subsec:Qscheme}
Let us  consider the {\em ARCH} model
\begin{equation}\label{eq:ARCHdef}
X_{k+1} = X_k +\vartheta_k(X_k) Z_{k+1},\;k=0:n-1,
\end{equation}
where the (Borel) functions $\vartheta_k:\R^d\to {\mathbb M}_{d,q}$ are locally bounded  and the r.v. $(Z_k)_{k=1:n}$  are square integrable, centered and mutually  independent (when   ${\rm Cov}(Z_k) = I_q$ for every $k=1:n$,  $(Z_k)_{k=1:n}$  is a white noise). 

\smallskip If the r.v. $Z_k$ all lie   in $L_{\R^q}^{\infty}(\P)$, then Assumption $({\cal G}_{G,Z})$ is satisfied since $G_k(x,z) = x+ \vartheta_k (x) z$. However such is not the case in many applications. Indeed, 
%
the {\em Euler scheme} with Brownian increments of a martingale Brownian diffusion with diffusion coefficient $\vartheta(t,x)$ is an {\em ARCH} model corresponding  to the choice $\vartheta_k(x)= \sqrt{\frac Tn}\vartheta(t_k, x)$ with a Gaussian  ${\cal N}(0;I_q)$-distributed white noise  $(Z_k)_{k=0:n}$ since
\begin{equation}\label{eq:Euler}
\bar X_{k+1} = \bar X_k +\vartheta(t_k, \bar X_k) \sqrt{\tfrac Tn } Z_{k+1}, \; k=0, \ldots,n-1
\end{equation}
where $t_0=0$, $t_k= \tfrac{kT}{n}$, $Z_k = \sqrt{\tfrac nT }(W_{t_k}-W_{t_{k-1}})\sim{\cal N}(0;I_q)$, $k=1:n$ ($W$ is a standard $q$-dimensional Brownian motion). However,  such a Gaussian  white noise makes impossible  {\em assumption~$({\cal G}_{G,Z})$  to hold true}.  






More generally, an {\em ARCH} model evolving according to~\eqref{eq:ARCHdef} with non-vanishing functions $\vartheta_k$ satisfies~\eqref{eq:HypoBounded} iff the noise $(Z_k)_{k=1:n}$ is compactly supported. In order to produce a finitely supported approximation scheme based on dual quantization, we need to perform a first approximation consisting in replacing the original white noise $(Z_k)_{k=1:n}$ by a {\em bounded} noise $(\breve{Z_k})_{k=1:n}$.  
 
 Then, as  second step we can dually quantize this new {\em ARCH} following~\eqref{eq:def0} and~\eqref{eq:def2}.
 
 
\medskip
 \noindent {\sc Step~1} ({\em Bounding the noise}). To bound the white noise $(Z_k)_{k=1:n}$ we map it into a sequence $(\breve{Z}_k)_{k:1:n}$ satisfying for every $k=1:n$,
 \begin{equation}\label{eq:breveZbound}
\breve{Z}_k = \varphi_{k}(Z_k) \; \mbox{in such a way that}\; \left\{\begin{array}{ll} (i)&\varphi_k: \R^d\to \R^d \mbox{ is Borel and bounded}\\
									 (ii)& \E\, \breve Z_k = 0,\\
                                                                              (iii) &\forall\, i, j \!\in \{1, \ldots, d\}, \; Z^i_k-\breve Z^i_k \perp_{L^2} \breve Z^j_k
                                                                            \end{array}\right.\; k=1:n.
\end{equation}
 
  Note that, although the random vectors $\breve Z_k$ are still independent by construction, the sequence  $(\breve Z_k)_{k=1:n}$ is not a white noise~--~except if $\breve Z_k = Z_k$ a.s.~--~due to $(iii)$ since $\| Z_k\|_2^2 = \|Z_k-\breve Z_k\|_2^2 +\|\breve Z_k\|_2^2$. In particular $\|\breve Z_k\|_2 \le \|Z_k\|_2$ with equality iff $Z_k = \breve Z_k$. A sequence satisfying~\eqref{eq:breveZbound} will be called a {\em quasi-white noise} in what follows and two canonical examples are given just below. 

 Then we define the {\em ARCH} model associated to $(\breve Z_k)_{k=1:n}$ (with the  diffusion coefficients $\vartheta_k$) by  
\begin{equation}\label{eq:ARCHTruncdef}
\breve X_{k+1} = \breve X_k +\vartheta_k(\breve X_k) \breve Z_{k+1} ,\; k=0:n-1,\quad \breve X_0= g_0(X_0)
\end{equation}
where $g_0:\R^d\to \R^d$ is a bounded Borel function (usually a Voronoi quantization of $X_0$ when $X_0$ is not deterministic). 

\smallskip It is clear that if ${\cal F}_k=\sigma(X_0,(Z_j)_{j=1:k})$ then
\begin{equation}\label{eq:breveXmart}
(\breve X_k)_{k=0:n}\; \mbox{  is an ${\cal F}_k$-martingale and an ${\cal F}_k$-Markov chain}.
\end{equation} 

\smallskip
\noindent  {\bf Examples of interest.} 
$\blacktriangleright$ {\em Truncated white noise.} Set 
\begin{equation}\label{eq:truncWN}
\breve Z_k = Z_k\mbox{\bf 1}_{\{Z_{k}\in A_{k}\}}, \; k=1,\ldots,n
\end{equation}
where $A_1,\ldots,A_n$ are compact sets of $\R^q$ such that 
\begin{equation}
   \E\, Z_{k}\mbox{\bf 1}_{\{Z_{k}\in A_{k}\}}=0,\; k=1,\ldots,n.\label{eqcentronc}
\end{equation}
Such Borel sets $A_k$ are easy to specify when the r.v. $Z_k$ have symmetric distributions  (invariant by multiplication by $-1$) since balls centered at $0$  (or any symmetric set)  are admissible. Notice that~\eqref{eq:breveZbound}$(iii)$ is satisfied whatever the choice of the sets $A_k$ whereas~\eqref{eqcentronc} is equivalent to \eqref{eq:breveZbound}$(ii)$.
 
 \smallskip
 \noindent $\blacktriangleright$ {\em Primal/Voronoi stationary quantization of the white noise.}  We replace the white noise by a quantization, usually a Voronoi (primal) one,  since the original white noise has no reason to be bounded. Then we set
\begin{equation}\label{eq:VorQuantZ}
 \breve Z_k = \widehat Z_k^{vor, \Gamma^Z_k}
\end{equation}
where $\Gamma^Z_k$ is a stationary primal/Voronoi  quantization grid with size $N^Z_{k}\ge 1$. Conditions~\eqref{eq:breveZbound}$(ii)$-$(iii)$  follow from the (primal) stationarity property~(see \eqref{eq:StatioVoro})
\begin{equation}
   \E\left[Z_k| \widehat Z_k^{vor, \Gamma^Z_k}\right]= \widehat Z_k^{vor, \Gamma^Z_k}\label{eq:VoroStatio}. 
\end{equation}

\smallskip
In both  settings, defining  $(\breve X_k)_{k=0:n}$ by~\eqref{eq:ARCHTruncdef}, one obtains an  approximation of $(X_k)_{k=0:n}$ which is both  non-decreasing for the convex order (as a martingale) and, as will be seen in Section~\ref{secpathconv}, is   dominated by the original {\em ARCH} model $(X_k)_{k=0:n}$ (under some additional convexity assumptions on the functions $\vartheta_k$ to be  made precise).

\medskip
\noindent{\sc STEP~2} ({\em Dual quantization}).
 We suppose we are given  a sequence of (finite) grids $\Gamma_k \subset \R^d$, $k=0,\ldots,n$  satisfying $\Gamma_0= g_0(\R^d)$ and the property $({\cal G}_{\Gamma, \breve Z})$, namely
 \[
(\mathcal{G}_{\Gamma, \breve Z})\;\;\; \forall\, x\!\in \Gamma_{k-1}, \quad x+\vartheta_{k-1}(x)\breve Z_k(\Omega)\subset {\rm conv}\big(\Gamma_k\big), \quad  k= 1:n.
\]
 Then we rely on the dual quantization procedure~\eqref{eq:def0} and~\eqref{eq:def2} to define recursively  the quantized approximation $(\widehat X_k)_{k=0:n}$  of $(\breve{X}_k)_{k=0:n}$ starting from $\widehat X_0= g_0(X_0)$ (usually a Voronoi quantization of $X_0$) supported by  a grid $\Gamma_0$ with $N_0$ points). Then 
\begin{equation}\label{eq:XZtildebreve}
\widetilde X_{k+1} = \widehat X_k+\vartheta_k(\widehat X_k)\breve Z_{k+1}\quad \mbox{ and }\quad \widehat X_{k+1}= {\rm Proj}^{del}_{\Gamma_{k+1}}(\widetilde X_{k+1}, U_{k+1}), \; k=0, \ldots,n-1
\end{equation}
where $(U_k)_{k=1:n}$ is a sequence of independent random variables uniformly distributed on $[0,1]$ independent from $(Z_k)_{k=1:n}$ and $X_0$ (hence $(Z_k, \breve Z_k)_{k=1:n}$ and $\widehat X_0$).

The algorithms performing the optimization of primal and dual quantization grids in general dimensions are respectively discussed in Appendices~\ref{App:A1} and~\ref{App:A2}. For a convenient implementation of the dual quantization of $\widehat X_k+\vartheta_k(\widehat X_k)\breve Z_{k+1}$, one needs to know the affine space spanned by the support of the distribution of this random vector. Let $\breve {\mathbb C}_{k+1}$ denote the covariance matrix of $\breve Z_{k+1}$. If for each $x\in\R^d$, the matrix $\vartheta_k(x)\breve{\mathbb C}_{k+1}\vartheta^*_k(x)$ is positive definite (which implies that $q\ge d$), then this affine space is $\R^d$.

\begin{Remark} Truncation of the noise is necessary to perform dual quantization of the chain. For numerical aspects on the way to discretize (or not) the truncated noise we refer to Section~\ref{sec:euler} where this problem is discussed in details in the case of the Euler scheme of a diffusion.
\end{Remark}

\subsection{Error estimates}
We equip the space $ \mathbb{M}_{d,q}$ with the operator norm  $\vertiii{B} = \sup_{|x|\le 1} |Bx|$ where $|\cdot|$ denotes the canonical Euclidean norm. Then for an $\mathbb{M}_{d,q}$-random variable $M$ we denote in short $\|M\|_2$ for $\big\| \,\vertiii{M}\,\big\|_2$. Then we will denote by  $[\vartheta]_{\rm Lip}$ the Lipschitz coefficient (if finite) of $\vartheta: (\R^d, |\cdot|)\to (\mathbb{M}_{d,q}, \vertiii{\cdot})$.  We will also make use of the Fr\"obenius norm ($\|B\|_{\rm Fr}= \sqrt{{\rm Tr}(BB^*)}$, $B\!\in \mathbb{M}_{d,q}$)  which satisfies $\vertiii{B}\le \|B\|_{\rm Fr} \le \sqrt{d\wedge q} \,\vertiii{B}$.

For a Lipschitz continuous function  $\vartheta:\R^d \to \mathbb M_{d,q}$, we define 
\begin{align*}
  c(\vartheta)= \sup_{x\in \R^d} \frac{\vertiii{\vartheta(x)}^2}{1+|x|^2}  \le 2\big(\vertiii{\vartheta(0)}^2 + [\vartheta]_{\rm Lip}^2\big)<+\infty\\
  c_{\rm Fr}(\vartheta)= \sup_{x\in \R^d} \frac{\|\vartheta(x)\|_{\rm Fr}^2}{1+|x|^2}  \le 2\big(\|\vartheta(0)\|_{\rm Fr}^2 + [\vartheta]_{\rm Fr, Lip}^2\big)<+\infty
\end{align*}
where $[\vartheta]_{\rm Fr, Lip}$  denotes the Lipschitz constants of $\vartheta$ with respect to the Fr\"obenius norm.

\smallskip
Thus, using that, if $\zeta \!\in L^2_{\R^q}(\P)$ is centered  with $L^2$-orthogonal components $\zeta^i$ and $B\in\mathbb M_{d,q}$, then $\E\, |B\zeta |^2 = \sum^d_{i=1}\sum^q_{j=1}B^2_{ij}\E\, (\zeta^j)^2$,
we straightforwardly derive the following inequality
\begin{align*}
\E\, |X_{k+1}|^2  = \E\, |X_k|^2 +\E\, |\vartheta_k(X_k) Z_{k+1}|^2&=  \E\, |X_k|^2 +\E\,|{\rm Tr}( \vartheta_k \vartheta_k^*)(X_k)|\\
&\le  \big(1+c_{\rm Fr}(\vartheta_k)\big) \E\, |X_k|^2  +c_{\rm Fr}(\vartheta_k).
\end{align*}
Standard induction then yields
\begin{equation} \label{eq:Bound0}
\|X_k\|_2^2 = \E\, |X_k|^2 \le \left[ \prod_{\ell=0}^{k-1} \big(1+c_{\rm Fr}(\vartheta_\ell)\big)\right] \Big(\E\, |X_0|^2+1\Big) -1
\end{equation}
so that
\begin{equation}
   \|\vartheta_k(X_k)\|_2^2\le c(\vartheta_k)^2\prod_{\ell=0}^{k-1} \big(1+c_{\rm Fr}(\vartheta_\ell)\big) \Big(\E\, |X_0|^2+1\Big).\label{eq:majotxk}
\end{equation}

We have the following strong $L^2$-error estimates where we can use this estimation of $\|\vartheta_k(X_k)\|_2^2$.

\begin{Thm} \label{thm:cvgce}Let $(U_k)_{k=0:n}$ be a sequence of i.i.d. random variables uniformly distributed over $[0,1]$ and let $(Z_k,\breve Z_k)_{k=1:n}$ be an independent sequence of independent square integrable $\R^{q+q}$-valued random vectors, independent of $X_0$  satisfying~\eqref{eq:breveZbound} and such that $(Z_k)_{k=1:n}$ is a white noise. Let $\widehat X_0=g_0(X_0)$ and let $(\Gamma_k)_{k=1:n}$ be grids satisfying  the consistency condition~$({\cal G}_{\Gamma,Z})$. 

\smallskip
\noindent $(a)$ {\em General case}. Then, for every $k\!\in \{0, \ldots,n\}$, 
\begin{align}
\hskip-0,15cm\| \widehat X_{k}-X_{k}\|_2^2&\le\left[ \prod_{i=1}^{k}\big(1+q[\vartheta_{i-1}]^2_{\rm Lip} \big)\right]\|\widehat X_{0}-X_{0}\|_2^2\nonumber \\
 &\quad +\sum_{\ell=1}^k\!\left[\prod_{i=\ell+1}^{k}\big(1+q[\vartheta_{i-1}]^2_{\rm Lip} \big) \right] \!  \times \!\bigg(\|\vartheta_{\ell-1}(X_{\ell-1})\|_2^2\big\|Z_{\ell}- \breve Z_{\ell}\big\|_2^2+\big\|\widehat X_{\ell}-\widetilde X_{\ell}\big\|_2^2\bigg).\label{eq:Upperfinala0}
\end{align}

Moreover, one has:  $\quad \displaystyle \Big\| \max_{k=0:n}|X_k-\breve X_k|\Big\|^2_2 \le 4 \big\| X_n-\breve X_n\big\|^2_2.$

\smallskip
\noindent $(b)$  {\em Quantized innovation.} Assume that $\widehat X_0=   {\rm Proj}^{vor}_{\Gamma_{0}}(X_0)$ is an optimal quadratic Voronoi quantization of $X_0$ at level $N_0$ and $\breve Z_k =  \widehat Z^{vor}_k = {\rm Proj}^{vor}_{\Gamma^Z_{k}}(Z_k)$, $k=1,\ldots,n$,   all are  optimal quadratic Voronoi quantizations of $Z_k$  at level $N^Z_k= {\rm card}(\Gamma^Z_k)$. Assume $X_0$, $Z_1, \ldots, Z_n \!\in L^{2+\eta}(\P)$. Then $X_k\!\in L^{2+\eta}(\P)$ for every $k\!\in \{0, \ldots,n\}$ and, for every $k=0, \ldots,n$,
\begin{align*}
  \|\widehat X_{k}&-X_{k}\|^2_2 \le (C^{vor}_{d,\eta})^2\prod_{i=1}^{k}\big(1+q[\vartheta_{i-1}]^2_{\rm Lip}\big) \frac{\sigma^2_{2+\eta}(X_0)}{N_0^{2/d}}
\\
 & \quad +\sum_{\ell=1}^k\left[\prod_{i=\ell+1}^k\hskip-0.25cm \big(1+q[\vartheta_{i-1}]^2_{\rm Lip} \big) \right] \left[ \|\vartheta_{\ell-1}(X_{\ell-1})\|_2^2(C^{vor}_{q,\eta})^2\frac{\sigma_{2+\eta}^2(Z_\ell)}{(N^Z_\ell)^{2/q}} + (\widetilde C^{del}_{d,\eta})^2 \frac{\sigma_{2+\eta}^2(\widetilde X_\ell)}{(N_\ell)^{2/d}}\right].
\end{align*}
\end{Thm}
\begin{Remark} \label{rk:diagq}
$\bullet$  If, furthermore, the $\breve Z_k$ have diagonal covariance matrices, 
then, for $A\!\in \mathbb{M}_{d,q}$,
\[
\E\, |A\breve Z_k|^2 = \sum_{i=1}^d\sum_{j=1}^q A^2_{ij}  \,\E\, \big(\breve Z^j_k\big)^2 \le\|A\|_{\rm Fr}^2\max _{j=1:q} \E\, \big(\breve Z^j_k\big)^2 \le \|A\|_{\rm Fr}^2.
\]
Moreover, for $i\neq j$ and $k=1,\ldots,n$,
$$\E[(\breve Z^i_k-Z^i_k)(\breve Z^j_k-Z^j_k)]=\E[(\breve Z^i_k-Z^i_k)\breve Z^j_k]+\E[\breve Z^i_k(\breve Z^j_k-Z^j_k)]+\E[Z^i_kZ^j_k-\breve Z^i_k\breve Z^j_k]
=0,$$
by~\eqref{eq:breveZbound} and since the covariance matrices of both $Z_k$ and $\breve Z_k$ are diagonal. Hence, for $A\!\in \mathbb{M}_{d,q}$,
\[
\E\, |A(Z_k-\breve Z_k)|^2 \le\|A\|_{\rm Fr}^2\max _{i=1:q} \E\, \big(Z_k^i-\breve Z^i_k\big)^2 .
\] 
Consequently \eqref{eq:Pyth} can be replaced by \begin{align*}
\nonumber\big \|(\vartheta_k(\widehat X_k)-\vartheta_k(X_k))\breve Z_{k+1}&+\vartheta_k(X_k)(\breve Z_{k+1}-Z_{k+1})\big\|_2^2\\
&\le  [\vartheta_{k}]^2_{\rm Fr,  Lip}\|\widehat X_k-X_k\|^2_2+\big\|\|\vartheta_k(X_k)\|_{\rm Fr}\big\|_2^2\max_{i=1:q} \E\,\big(Z^i_{k+1}- \breve Z^i_{k+1}\big)^2.
\end{align*} 
Hence~the estimations in Theorem \ref{thm:cvgce} still hold with $q[\vartheta_{i-1}]^2_{\rm Lip}$ and $\|\vartheta_{\ell-1}(X_{\ell-1})\|_2^2\|Z_{\ell}-\breve Z_{\ell}\|_2^2$  replaced by $[\vartheta_{i-1}]^2_{\rm Fr,  Lip}$,  $\big\|\|\vartheta_{\ell-1}(X_{\ell-1})\|_{\rm Fr}\big\|_2^2\max_{i=1:q}\E\,\big(Z^i_{\ell}- \breve Z^i_{\ell}\big)^2$ respectively.

\smallskip
\noindent $\bullet$ Let us assume that the vectors $Z_k$ have {\em independent coordinates}. Then, the diagonal covariance matrix condition can be achieved by choosing 
$$
\breve{Z}^i_k = \varphi_{k,i}(Z^i_k)
$$
 in such a way that~\eqref{eq:breveZbound}$(ii)$ is satisfied  and~\eqref{eq:breveZbound}$(iii)$ holds for $i=j$ (for every $i$). Then, for $i\neq j$,~\eqref{eq:breveZbound}$(iii)$ follows from the independence property.  

\smallskip
As for truncation, this can be done by considering product sets of the form  $A_k =\prod_{i=1}^qA^i_k$ such that $\E\, Z^i_k \mbox{\bf 1}_{\{Z^i_k \in A^i_k\}}= 0$, $i=1:q$. 

As for the quantization based approach, first note that optimal Voronoi quantization usually does not satisfy this property. But this can be achieved by calling upon {\em product quantization}: one  considers product grids $\Gamma^Z_k= \prod_{1\le i \le q}\!\Gamma^Z_{k,i}$ where $\Gamma^Z_{k,i}$ is a  stationary  Voronoi grid of the $i$-th marginal $Z^i_k$ so that 
$$
\E \big(Z^i_k |(\widehat Z_k^{j,vor})_{j=1:q}\big)=\E\,\big(  Z^i_k\, |\,\widehat Z^{i,vor}_k\big) = \widehat Z^i_k,\quad i=1,\cdots,q,
$$  
where $\widehat Z^{i,vor}_k = {\rm Proj}_{\Gamma^Z_{k,i}}^{vor}(Z^i_k)$. One easily derives from this componentwise stationarity property that~\eqref{eq:breveZbound} holds true for $\breve Z_k= \big(\widehat Z^{i,vor}_k\big)_{i=1:q} = \widehat Z^{\Gamma_{_{\!Z}} ,vor}_k$. Finally note that, when these marginal quantizations $\widehat Z^{i,vor}_k$ are $L^2$-optimal, then  $\widehat Z^{\Gamma_{_{\!Z}} ,vor}_k$ is rate optimal and satisfies the universal non-asympotic upper-bound provided by the so-called Pierce Lemma (see the remark following Theorem ~\ref{thm:zadoretpierce1} in Appendix~\ref{App:A1}).
\end{Remark}

\begin{Remark}\label{remconvgen}Let us consider the approximation \eqref{eq:def0}-\eqref{eq:def2} of the more general model \eqref{eq:Markovchain}.
   Under the following $L^2$-Lipschitz assumption on the function $G_k$ and the r.v. $Z_{k+1}$,
\begin{equation}\label{eq:LipG_k}
\forall\, k=0, \ldots,n-1,\\; \forall\,x,y\!\in \R^d,\; \big\|G_k(x, Z_{k+1})-G_k(y,Z_{k+1}) \big \|_2 \le [G_k]_{\rm Lip}|x-y|
\end{equation}
with  $\max_{k=0,\ldots,n-1}[G_k]_{\rm Lip}<+\infty$, we can check in the same way that for every $k\!\in \{1,\ldots,n\}$
\[
\big\| \widehat X_{k}-X_{k}\big\|_2\le \left(\sum_{\ell=0}^k [G_{\ell:k}]^2_{\rm Lip}\big\| \widehat X_{\ell}-\widetilde X_\ell \big\|^2_2\right)^{1/2},
\]
 with the convention $\widetilde X_0= X_0$ and where, for $0\le \ell\le k$, $\displaystyle 
[G_{\ell:k}]_{\rm Lip} =\prod_{i=\ell+1}^k[G_i]_{\rm Lip}$ (with $\prod_{\varnothing} =1$).
\end{Remark}

\noindent {\bf Proof.} 
%
Both $(X_k)_{k=0:n}$ and    $(\widehat X_k)_{k=0:n}$ are $\sigma\big(X_0,(Z_{\ell} , U_{\ell})_{\ell  = 1: k} \big)$-martingales so that, on the one hand,  $\widehat X_k\le_{cvx}\widehat X_{k+1}$ for all $k=0, \ldots,n-1$. On the other hand, their difference is also a martingale and  one derives from the decomposition 
\begin{align*}
   \widehat X_{k+1}&-X_{k+1} =\widehat X_k-X_k\\
   &\quad + {\rm Proj}^{del}_{\Gamma_{k+1}}\big(\widetilde X_{k+1},U_{k+1}\big)-\widetilde X_{k+1}+\big(\vartheta_k(\widehat X_k)-\vartheta_k(X_k)\big)\breve Z_{k+1}+\vartheta_k(X_k)(\breve Z_{k+1}-Z_{k+1})
        \end{align*}
        that 
        \begin{align*}
  \| \widehat X_{k+1}&-X_{k+1}\|_2^2 = \|\widehat X_k-X_k\|_2^2\\
   &\quad + \|{\rm Proj}^{del}_{\Gamma_{k+1}}\big(\widetilde X_{k+1},U_{k+1}\big)-\widetilde X_{k+1}+\big(\vartheta_k(\widehat X_k)-\vartheta_k(X_k)\big)\breve Z_{k+1}+\vartheta_k(X_k)(\breve Z_{k+1}-Z_{k+1})\|_2^2.
         \end{align*}
        As the random variable $U_{k+1}$  is independent  of $(\widetilde X_{k+1}, \widehat X_k, X_k,  Z_{k+1}, \breve Z_{k+1})$, it follows from the dual stationarity property that 
 \begin{align}
\nonumber    \|\widehat X_{k+1}-X_{k+1}\|_2^2& =\|\widehat X_{k}-X_{k}\|_2^2+\big\|{\rm Proj}^{del}_{\Gamma_{k+1}}\big(\widetilde X_{k+1},U_{k+1}\big)-\widetilde X_{k+1}\big\|_2^2\\
   \label{eq:decapprox1}&\quad+\big\|(\vartheta_k(\widehat X_k)-\vartheta_k(X_k))\breve Z_{k+1}+\vartheta_k(X_k)(\breve Z_{k+1}-Z_{k+1})\big\|_2^2.
 \end{align}
Moreover, by~\eqref{eq:breveZbound}$(iii)$ and the independence of $(Z_{k+1}, \breve Z_{k+1})$ and $(X_k, \widehat X_k)$, one has
\begin{align}
\nonumber\big \|(\vartheta_k(\widehat X_k)-\vartheta_k(X_k))\breve Z_{k+1}&+\vartheta_k(X_k)(\breve Z_{k+1}-Z_{k+1})\big\|_2^2\\
\nonumber & = \|(\vartheta_k(\widehat X_k)-\vartheta_k(X_k))\breve Z_{k+1}\|_2^2+\|\vartheta_k(X_k)(\breve Z_{k+1}-Z_{k+1})\|_2^2\\
\label{eq:Pyth}&\le  \|(\vartheta_k(\widehat X_k)-\vartheta_k(X_k))\|^2_2\|\breve Z_{k+1}\|^2_2+\|  \vartheta_k(X_{k})\|_2^2 \|Z_{k+1}- \breve Z_{k+1}\|_2^2.
\end{align}

%
  It follows from~\eqref{eq:decapprox1},~\eqref{eq:Pyth}, the Lipschitz property of the functions  $\vartheta_k$ and  the inequality $\| \breve Z_{k+1}\|_2^2  \le  \| Z_{k+1}\|_2^2=q$ deduced from Condition~\eqref{eq:breveZbound}$(iii)$ that
  \begin{align*}
  \|\widehat X_{k+1}-X_{k+1}\|_2^2
 & \le \|\widehat X_{k}-X_{k}\|_2^2 \big(1+q[\vartheta_k]^2_{\rm Lip}\big) +  \|  \vartheta_k(X_{k})\|_2^2 \|Z_{k+1}- \breve Z_{k+1}\|_2^2
\\
&\quad +\|\widehat X_{k+1}-\widetilde X_{k+1}\|_2^2.
  \end{align*}
The discrete time Gronwall's lemma yields for every $k=0, \ldots,n$,
\begin{align}
\| \widehat X_{k}-X_{k}\|_2^2&\le \left[\prod_{i=1}^{k}\big(1+q[\vartheta_{i-1}]^2_{\rm Lip} \big)\right]\|\widehat X_{0}-X_{0}\|_2^2\nonumber \\
 &\hskip 0,15cm +\sum_{\ell=1}^k\left[\prod_{i=\ell+1}^{k}\big(1+q[\vartheta_{i-1}]^2_{\rm Lip} \big) \right] \! \times\!\bigg(\|\vartheta_{\ell-1}(X_{\ell-1})\|_2^2\big\|Z_{\ell}- \breve Z_{\ell}\big\|_2^2+\big\|\widehat X_{\ell}-\widetilde X_{\ell}\big\|_2^2\bigg).\label{eq:Upperfinala}
\end{align}

\smallskip
\noindent $(b)$ Optimality of the quantizations $\breve Z_{k}$ imply that  the $\breve Z_{k}$ are Voronoi stationary which makes  $(Z_k, \breve Z_k)$ a martingale coupling hence satisfying~\eqref{eq:breveZbound}  for every $k=1,\ldots,n$. As $X_0\in L^{2+\eta}_{\R^d}(\P)$ and $Z_k \!\in L^{2+\eta}_{\R^q}(\P)$ for $k=1,\ldots,n$,  the Voronoi (primal) non-asymptotic version of Zador's Theorem (see Theorem~\ref{thm:zadoretpierce1}$(b)$ in Appendix~\ref{App:A1}) implies that
\[
\big\|X_{0} - \widehat X_{0}\big\|_2 \le \widetilde C^{vor}_{d,\eta} \sigma_{2+\eta}(X_0)N_0^{-1/d}\;\mbox{ and }\;\big \|Z_{k} - \breve Z_{k}\big\|_2 \le \widetilde C^{vor}_{2,q,\eta} \sigma_{2+\eta}(Z_k)(N^Z_k)^{-1/q},\quad k=1,\ldots,n,
\]
where $\widetilde C^{vor}_{2,q,\eta}$ is a positive real constant only depending on the dimension $q$ and $\eta>0$. 
Moreover, for every $k=1,\ldots,n$, the random variables $\widetilde X_k$ are compactly supported. Hence, owing to the dual form of Zador's Theorem (see Appendix~\ref{App:A2}, Theorem~\ref{thm:zadoretpierce2}$(b)$), there exists a real constant $\widetilde C_{2,d, \eta}^{del}\!\in (0, +\infty)$ such that, for every $k=1,\ldots,n$, 
\[
\big\|\widehat X_{k} -\widetilde X_{k}\big\|_2\le \widetilde C^{del}_{2,d,\eta}\sigma_{2+\eta}(\widetilde X_{k})N_k^{-1/d}.
\]

Plugging these bounds  into~\eqref{eq:Upperfinala} completes the proof. $\hfill\cqfd$

\medskip It should be noted  that ``skipping'' the dual quantization step in the above proof and using \eqref{eq:majotxk} also provides the following convergence rate for $(\breve X)_{k=0:n}$ defined by \eqref{eq:ARCHTruncdef} toward the original {\em ARCH} $(X_k)_{k=0:n}$.

 \begin{Pro} \label{prop:3.1}Assume that all the functions $\vartheta_k:\R^d\to \mathbb{M}_{d,q}$, $k=0:n-1$, are Lipschitz continuous and that $X_0 \!\in L^{2}_{\R^d}(\P)$ ($e.g.$ because $X_0=x_0\!\in \R^d$). Then, for every $k=0:n$, 
\begin{align}
\label{eq:Bound1}
 \hskip-0,25cm \|X_k-\breve X_k\|_2^2 & \le    \| X_0-\breve X_0\|_2^2 \prod_{i=1}^{k} (1+ q[\vartheta_{i-1}]_{\rm lip}^2)\! +\!  \big(1+\|X_0\|^2_2\big) \!\!
 \prod_{i=1}^{k}\!\!\big(1+C(\vartheta_{i-1})\big) \!\sum_{\ell=1}^{k}c(\vartheta_{\ell-1}) \|Z_{\ell}-\breve Z_{\ell}\|_2^2,
\end{align}
where $C(\vartheta) = \big(q[\vartheta]_{\rm Lip}^2\big)\vee c_{\rm Fr}(\vartheta)$. Moreover, $
\displaystyle \Big\| \max_{k=0:n}|X_k-\breve X_k|\Big\|^2_2 \le 4 \big\| X_n-\breve X_n\big\|^2_2$.
\end{Pro}
%

\begin{Remark} \label{rk:diagq}
If, furthermore, the $\breve Z_k$ have diagonal covariance matrices, 
then~\eqref{eq:Bound1} holds with $C(\vartheta) = [\vartheta]^2_{\rm Fr,  Lip}\vee c_{\rm Fr}(\vartheta)$ and with $q[\vartheta_{i-1}]^2_{\rm Lip}$ , $c(\vartheta_{\ell-1}) \|Z_{\ell}-\breve Z_{\ell}\|_2^2$ respectively replaced by $[\vartheta_{i-1}]^2_{\rm Fr,  Lip}$,  $c(\vartheta_{\ell-1})_{\rm Fr} \max_{i=1:q}\E\big(Z^i_{\ell}-\breve Z^i_{\ell}\big)^2$.
\end{Remark}
\section{Convex {\em ARCH} models: }\label{sec:convex}
In this section, we are first going to give conditions on the functions $(\vartheta_k)_{k=0:n-1}$ and the noise $(Z_k)_{k=1:n}$ ensuring that $(\breve X_k)_{k=0:n}\le_{cvx}(X_k)_{k=0:n}$. Since each dual quantization step is convex order increasing, $(\widehat X_k)_{k=0:n}$ is, in general, not comparable to the original ARCH $(X_k)_{k=0:n}$. Nevertheless, the fact that  truncation and  dual quantization have opposite effects in terms of convex order is not so bad for numerical purposes: the errors coming from these two approximations should, at least partially, compensate.
Moreover, we will also check that under the same conditions on the functions $(\vartheta_k)_{k=0:n-1}$ and the truncated noise $(\breve Z_k)_{k=1:n}$, $(\breve X_k)_{k=0:n}\le_{cvx}(\widehat X_k)_{k=0:n}$.

\subsection{Convex domination of $(\breve X_k)_{k=0:n}$ by $(X_k)_{k=0:n}$}\label{secpathconv}
When the functions $\vartheta_k$ are convex in an appropriate sense and  the variables $Z_{k+1}$ have radial distributions, then the {\em ARCH} model~\eqref{eq:ARCHdef} dominates all its approximations  $(\breve X_k)_{k=0:n}$ with  truncated or quantized white noise. This is summed in the next theorem, stated right after the definition of convexity we need when dealing with matrix-valued functions.

\begin{Dfn} A function $\vartheta : \R^d \to \mathbb{M}_{d,q}(\R)$ is convex if 
  \begin{align}
\forall\, x,\, y\!\in \R^d,\forall &\alpha\in[0,1],\;\exists \, O = O_{x,y,\alpha}\in{\mathbb M}_{q,q}\mbox{ orthogonal such that }\notag\\&\vartheta \vartheta^*\big(\alpha x +(1-\alpha) y\big)\le \big(\alpha \vartheta(x) +(1-\alpha) \vartheta(y)  O\big)\big (\alpha \vartheta(x) +(1-\alpha) \vartheta(y)O \big)^*.
\label{eq:convmat} 
\end{align}
\end{Dfn}

Note that when $d=q=1$, the left-hand side is equal to $|\vartheta(\alpha x +(1-\alpha) y)|^2$ while the right-hand side is maximal and equal to $\big(\alpha |\vartheta(x)| +(1-\alpha) |\vartheta(y)|)^2$ when $O$ is equal to the sign of $\vartheta(x)\vartheta(y)$. Therefore the above convexity assumption is equivalent to the convexity of $|\vartheta|$.

\begin{Thm}[Domination]\label{propdomgausd} Let $(X_k)_{k=0:n}$ be an $\R^d$-valued {\em ARCH} model defined by~\eqref{eq:ARCHdef} where the $\R^q$-valued white noise $(Z_k)_{k=1:n}$ is a sequence of 
  $\R^q$-valued r.v. with radial distributions, the initial random vector $X_0$ is integrable and the ${\mathbb M}_{d,q}$-valued functions $\vartheta_k$, $k=0,\ldots,n-1$, are convex in the sense of~\eqref{eq:convmat} with linear growth. Assume that $g_0(X_0)\le_{cvx}X_0$.

\smallskip
\noindent $(a)$ {\em Truncation.} Let  $(\breve Z_k=Z^{A_k}_k)_{k=1:n}$ with $(A_k)_{k=1:n}$
an $n$-tuple of Borel sets satisfying 
$ \E\, Z_k\mbox{\bf 1}_{\{Z_k \in A_k\}}=0$, $k=1,\ldots,n$ and  $(\breve X^A_k)_{k=0:n}$  be  the induced approximating {\em ARCH} process defined by~\eqref{eq:ARCHTruncdef}. Then
\[
\breve X^A_{0:n} \le_{cvx}
X_{0:n}.
\]  
\noindent $(b)$ {\em Quantization.} Let $(\breve Z_k= \widehat Z^{vor}_{k})_{k=1:n}$ be a stationary (Voronoi) quantized approximation of the white noise $Z_{1:n}$ and $\breve X_{0:n}$ be the induced approximating {\em ARCH} process defined by~\eqref{eq:ARCHTruncdef}. Then 
\[
\breve X_{0:n} \le_{cvx} X_{0:n}.
\]

\noindent $(c)$ {\em Scalar setting.} When $d=q=1$, if, mutatis mutandis,  $|\vartheta_k|$ is convex for every $k=0,\ldots,n-1$, then the conclusions of $(a)$ and~$(b)$ hold true. 
\end{Thm}

\smallskip 
First, we check that that truncation and quantization induce a convex domination. For Voronoi quantization, this is clear by the stationary property since,  if $\breve{Z}_k = \widehat Z^{vor}_k$, then $\breve{Z}_k\le_{cvx} Z_k$, $k=1,\ldots,n$.

\smallskip
The truncation case is solved by the  following lemma.
\begin{Lem}[Truncation]\label{lem:truncZ} Let $Z\!\in L_{\R^q}^1(\Omega,{\cal A}, \P)$ be a centered random vector. For any Borel subset $A$ of $\R^q$, let $ Z^A= Z\mbox{\bf 1}_{\{Z\in A\}}$. If $\E\, Z^A=0$, then 
\[
Z^A\le _{cvx} Z.
\] 
\end{Lem}

\noindent {\bf Proof.} One may restrict to convex functions $\varphi:\R^q\to \R$ with linear growth (see $e.g.$ Lemma A.1 in~\cite{AlJo} or Remark 1.1 p.2 in~\cite{HiPrRoYo}) for which $\varphi(Z)\!\in L^1$. We may assume w.l.g. that $\P(Z\notin A)>0$ (otherwise the result is trivial). Then
\begin{align*}
\E\, \varphi(Z) -\E\, \varphi(Z^A)& = \E\, \varphi(Z)\mbox{\bf 1}_{\{Z\notin A\}}-\varphi(0)\P(Z\notin A)\\
&= \P(Z\notin A) \Big(\E\,\big(\varphi(Z)\,|\, Z\notin A\big)-\varphi(0)\Big)\ge 0
\end{align*}
owing to Jensen's Inequality and $\E\, (Z\, |\, Z\!\notin A) = 0$.~$\hfill\cqfd$

Since, under the assumptions of Theorem \ref{propdomgausd}, both the quantized and the truncated approximations of the noise are smaller then this noise in the convex order, the theorem is a consequence of the following proposition applied with $(Z'_{k+1},\vartheta'_k)_{k=0:n-1}=(\breve Z_{k+1},\vartheta_k)_{k=0:n-1}$.
\begin{Pro}[Convex order:  from the noise to the {\em ARCH}]\label{Proconvgen} Let $(Z_k)_{k=1:n}$ and $(Z'_k)_{k=1:n}$ be two sequences of $\R^q$-valued independent integrable and centered random vectors. Let $(\vartheta_k)_{k=0:n-1}$ and $(\vartheta'_k)_{k=0:n-1}$ be two sequences of ${\mathbb M}_{d,q}$-valued functions with linear growth defined on $\R^d$ such that: $\forall x\in\R^d$, $\vartheta'_k(\vartheta_k')^*(x)\le \vartheta_k\vartheta_k^*(x)$ for $k=0, \ldots,n-1$ and, either for every $k=0, \ldots,n-1$,   $\vartheta_k$ is convex in the sense of~\eqref{eq:convmat}
and 
  the r.v. $Z_{k+1}$ has a radial distribution ({\bf we say that the assumption is satisfied by $(Z_{k+1},\vartheta_k,\vartheta'_k)_{k=0:n-1}$}) or 
  for each $k=0, \ldots,n-1$,  $\vartheta_k'$  is convex in the sense of~\eqref{eq:convmat} and $Z'_{k+1}$ has a radial distribution ({\bf we say that the assumption is satisfied by $(Z'_{k+1},\vartheta'_k,\vartheta_k)_{k=0:n-1}$}).
  
Let $X_0$ and $X'_0$ be integrable $\R^d$-valued r.v. independent of $(Z_k)_{k=1:n}$ and $(Z'_k)_{k=1:n}$ respectively.
 Denote by $(X_k)_{k=0:n}$ and $(X'_k)_{k=0:n}$ the two {\em ARCH} models respectively defined by~\eqref{eq:ARCHdef} and by $X'_{k+1}=X'_{k}+\vartheta'_k(X'_{k})Z'_{k+1}$ for $k=0, \ldots,n-1$.

\smallskip If $X'_0\le_{cvx} X_0$ and $Z'_k\le_{cvx} Z_k$ for every $k=1,\ldots,n$, then 
$$
(X'_k)_{k=0:n} \le_{cvx} (X_k)_{k=0:n}.
$$
\end{Pro}

We will see that, when $d=q$, this result is in fact a particular case of a more general result established at the end of this section where  we will consider two sequences $(Z_k)_{k=1:n}$ and $(Z'_k)_{k=1:n}$   of $\R^q$-valued independent integrable and centered random vectors such that $Z'_k \le_{cvx}Z_k$, $k=1,\ldots,n$, with possible non-radial distributions but only sharing some various lighter symmetry properties.  This is the reason why, on the way to the proof and in particular in the next lemma, we will deal with convex order for random vectors with more general distributions than radial ones. 

\begin{Lem}\label{pro:Cvxdq}Let $Z$ be an integrable (centered) $ \R^q$-valued r.v. For $i=1:q$,  denote by $Z_{-i}$ the subvector obtained by removing  the  $i$-th coordinate $Z_i$ from $Z$.
  \begin{description}
  \item[$(i)$ Vanishing conditional expectations assumption.] If for $i=1:q$, $\E[Z_i|Z_{-i}]=0$ a.s. and  $0\le \lambda_i\le \ell_i$ , then $${\rm Diag}(\lambda_{1},\hdots,\lambda_q) Z\le_{cvx}{\rm Diag}(\ell_{1},\hdots,\ell_q) Z$$ where ${\rm Diag}(\lambda_{1},\hdots,\lambda_q)\in{\mathbb M}_{q,q}$ denotes the diagonal matrix with diagonal elements $\lambda_1,\hdots\lambda_q$,
\item[$(ii)$ Symmetric conditional laws assumption.] If for each $i=1:q$, the conditional laws of $Z_i$ and $-Z_i$ given $Z_{-i}$ coincide a.s. and  $|\lambda_i|\le |\ell_i|$, then $${\rm Diag}(\lambda_{1},\hdots,\lambda_q) Z\le_{cvx}{\rm Diag}(\ell_{1},\hdots,\ell_q) Z.$$
  \item[$(iii)$ Radial distribution assumption.] If $A$, $B\!\in {\mathbb M}_{d,q}$ and $Z$ has a radial distribution  (for each orthogonal matrix $O\in{\mathbb M}_{q,q}$, $OZ$ has the same distribution as $Z$), then $AA^*\le BB^*\Rightarrow AZ\le_{cvx}BZ$. If moreover $\E|Z|^2\in (0,+\infty)$, then the converse implication holds.
  \end{description}

\end{Lem}
\begin{Remark} $\bullet$ 
  The radial distribution assumption implies the symmetric conditional law assumption (choose the orthogonal transformation which only inverts the sign of the $i$-th coordinate). This implies the vanishing conditional expectation assumption. On the other hand, the assumptions on the matrices multiplying $Z$ get weaker from $(i)$ to $(iii)$. 

\smallskip
\noindent $\bullet$  When $Z$ follows the (radial) distribution ${\mathcal N}(0;I_q)$ and $AA^*\le BB^*$ then for $\zeta\sim{\mathcal N}(0;BB^*-AA^*)$ independent of $Z$, $\E[AZ+\zeta|Z]=AZ$, so that $AZ\le_{cvx}AZ+\zeta$ and $AZ+\zeta\sim{\mathcal N}(0;BB^*)$ so that $AZ+\zeta$ has the same distribution as $BZ$. Hence $AZ\le_{cvx}BZ$. This is a simple alternative argument to that in~\cite{FadPag} which has inspired the generalization to any radial distribution below.
\end{Remark}
\noindent{\bf Proof.}$(i)$ 
For $X$ an integrable and centered random variable, and $\psi:\R\to\R$ a convex function with linear growth, the function $\R\ni u\mapsto \E\, \psi(uX)$ is clearly convex 	and attains its minimum at $u=0$ owing to Jensen's inequality. Hence it is non-decreasing on $\R_+$ and non-increasing on $\R_-$.
Now for $\varphi:\R^d\to\R$ convex with linear growth, repeatedly using the monotonicity property on $\R_+$, one obtains
\begin{align*}
   \E\,\varphi({\rm Diag}(\lambda_1,\hdots,\lambda_q)Z)&=\E\,\E\big((\varphi({\rm Diag}(\lambda_1,\hdots,\lambda_q)Z)|Z_{-1}\big)\le \E\,\E\big(\varphi({\rm Diag}(\ell_1,\lambda_2,\hdots,\lambda_q)Z)|Z_{-1}\big)\\&= \E\,\E(\varphi({\rm Diag}(\ell_1,\lambda_2,\hdots,\lambda_q)Z)|Z_{-2})\\
   & \le \E\,\E(\varphi({\rm Diag}(\ell_1,\ell_2,\lambda_3,\hdots,\lambda_q)Z)|Z_{-2})\le \hdots\le \E\,\varphi({\rm Diag}(\ell_{1},\hdots,\ell_q)Z).
\end{align*}
By Lemma A.1 in~\cite{AlJo} (see also Remark 1.1 p.2 in~\cite{HiPrRoYo}), 
one concludes that ${\rm Diag}(\lambda_{1},\hdots,\lambda_q)Z\le_{cvx}{\rm Diag}(\ell_{1},\hdots,\ell_q)Z$.

\smallskip
\noindent $(ii)$ Since, under the assumption, ${\rm Diag}(\lambda_{1},\hdots,\lambda_q)Z$ and ${\rm Diag}(\ell_{1},\hdots,\ell_q)Z$ respectively have the same distributions as ${\rm Diag}(|\lambda_{1}|,\hdots,|\lambda_q|)Z$ and ${\rm Diag}(|\ell_{1}|,\hdots,|\ell_q|)Z$, the conclusion follows from~$(i)$.

\smallskip
\noindent $(iii)$ {\bf Step 1}. For $C\in {\mathbb M}_{q,q}$, the singular value decomposition of $C$ writes $C=ODV$ for matrices $O,D,V\in {\mathbb M}_{q,q}$ with $O,V$ orthogonal and $D$ diagonal with nonnegative diagonal elements. One has $\sqrt{CC^*}=ODO^*$ and if $Z$ has a radial distribution, for any measurable and bounded function $\varphi:\R^q\to\R$, $\E\,\varphi(CZ)=\E\,\varphi(ODVZ)=\E\,\varphi(ODO^*Z)=\E\,\varphi(\sqrt{CC^*}Z)$ so that $CZ$ and $\sqrt{CC^*}Z$ share the same distribution.

\smallskip
\noindent{\bf Step 2}. Let us now assume that the $\R^q$-valued r.v. $Z$ has a radial distribution, $AA^*\le BB^*$ and $d=q$. We set $B_\varepsilon=\sqrt{BB^*+\varepsilon I_q}$. We have $B_\varepsilon^{-1}AA^*(B^{-1}_\varepsilon)^*\le I_q$ for $\varepsilon>0$. One deduces that $\sqrt{B_\varepsilon^{-1}AA^*(B^{-1}_\varepsilon)^*}=ODO^*$ for matrices $O,D\in {\mathbb M}_{q,q}$ with $O$ orthogonal and $D$ diagonal with diagonal elements belonging to $[0,1]$. For $\varphi :\R^d\to\R$ convex with linear growth, the function $\psi(x)=\varphi(B_\varepsilon O x)$ is convex with linear growth and
\begin{align*}
   \E\,\varphi(AZ)=\E\,\psi(O^*B_\varepsilon^{-1}AZ)=\E\,\psi(O^*ODO^*Z)=\E\,\psi(DZ)\le \E\,\psi(Z)=\E\,\varphi(B_\varepsilon OZ)=\E\,\varphi(B_\varepsilon Z),
\end{align*}
where we used the definition of $\psi$ for the first and fourth equalities, Step 1 for the second equality, the radial property of the distribution of $Z$ for the third and fifth equalities and $(i)$ for the inequality.
One has $\lim_{\varepsilon\to 0}B_\varepsilon=\sqrt{BB^*}$, so that by Lebesgue's theorem and Step 1, $\lim_{\varepsilon \to 0}\E\,\varphi(B_\varepsilon Z)=\E\,\varphi(\sqrt{BB^*}Z)=\E\,\varphi(BZ)$. We deduce that $\E\,\varphi(AZ)\le \E\,\varphi(BZ)$ so that $AZ\le_{cvx}BZ$.

\smallskip
\noindent{\bf Step 3}. Let us now assume that $Z$ has a radial distribution, $AA^*\le BB^*$ and $d<q$. Let $\widetilde A,\widetilde B\in {\mathbb M}_{q,q}$ be defined by 
$$
(\widetilde A_{ij},\widetilde B_{ij})=\begin{cases}(A_{ij},B_{ij})\mbox{ for }i=1:d,j=1:q\\
  (0,0)\mbox{ for }i=d+1:q,j=1:q\end{cases}.
  $$
We have $\widetilde A\widetilde A^*\le \widetilde B\widetilde B^*$, so that, by Step 2, $\widetilde A Z\le_{cvx}\widetilde BZ$. For $M\in{\mathbb M}_{d,q}$ with non-zero coefficients $M_{ii}=1$ for $i=1:d$, we have $AZ=M\widetilde A Z$ and $BZ=M\widetilde B Z$. Since for any convex function $\varphi:\R^d\to\R$, $\R^q\ni x\mapsto\varphi(Mx)$ is convex as the composition of a convex function with a linear function, we conclude that $AZ\le_{cvx}BZ$.

\smallskip
\noindent{\bf Step 4}. Let us finally assume that $Z$ has a radial distribution, $AA^*\le BB^*$ and $d>q$. 
We have $\mbox{Ker} B^*\subset \mbox{Ker}A^*$ so that $\mbox{Im} A\subset \mbox{Im} B$. Let $O\in{\mathbb M}_{d,q}$ be a matrix with orthogonal columns with norm one such that the first $\mbox{dim Im} B$ (we have $\mbox{dim Im} B\le q$) columns form an orthonormal basis of $\mbox{Im} B$. Then $B=OO^*B$ and $A=OO^*A$, $O^*BB^*O\ge O^*AA^*O$ and, by Step 2, $O^*AZ\le_{cvx}O^*BZ$. Since for any convex function $\varphi:\R^d\to\R$, $\R^q\ni x\mapsto\varphi(Ox)$ is convex as the composition of a convex function with a linear function, we conclude that $AZ=OO^*AZ\le_{cvx}OO^*BZ=BZ$.

\smallskip
\noindent{\bf Step 5}. Let us suppose that $Z$ is square integrable with a radial distribution. Then $\E(Z_iZ_j)=\mbox{\bf 1}_{\{i=j\}}\frac{\E|Z|^2}{q}$. If $A,B\in{\mathbb M}_{d,q}$ are such that $AZ\le_{cvx}BZ$, then, for $u\in\R^d$, the choice of the convex function $\varphi: x\in\R^d\mapsto (u^*x)^2$ in the inequality defining the convex order yields $u^*AA^*u\frac{\E|Z|^2}{q}\le u^*BB^*u\frac{\E|Z|^2}{q}$.~$\hfill\cqfd$

\medskip
\noindent {\bf Proof of Proposition \ref{Proconvgen}.} By the linear growth of the coefficients $\vartheta_k$, the integrability of the initial conditions and the noises and the independence structure, one easily checks by forward induction that $X_{0:k}$ and $X'_{0:k}$ are integrable for every $k=0, \ldots,n$. According to Lemma A.1 in~\cite{AlJo}, it is enough to prove that $\E\,\Phi_n(X_{0:n})\le \E\,\Phi_n(X'_{0:n})$ for $\Phi_n: (\R^d)^{n+1}\to \R$ convex with linear growth.

We proceed by successive backward  inductions. We define  the functions $\Phi_k: (\R^d)^{k+1}\to \R$, $k=0, \ldots,n-1$,  by backward induction as follows:
\[
\Phi_k(x_{0:k}) = \Psi_{k}\big(x_{0:k},\vartheta_k(x_k)\big), \; k=0, \ldots,n-1
\]
where, for every $(x_{0:k},u)\!\in (\R^d)^{k+1}\times {\mathbb M}_{d,q}$,
\[
\Psi_k(x_{0:k}, u) = \E\,\Phi_{k+1} \big(x_{0:k}, x_k+uZ_{k+1}\big), \; k=0, \ldots,n-1.
\]
By backward induction, using the integrability of the random variables $Z_{k+1}$ and the linear growth of the function $\vartheta_k(x_k)$, one easily checks that the functions $\Phi_k$ and $\Psi_k$ all have linear growth and in particular that the expectation in the definition of $\Psi_k$ makes sense.

Starting from $\Phi'_n=\Phi_n$, we define the functions $\Phi'_k$, $\Psi'_k$, $k=0, \ldots,n-1$ likewise using the  sequence $(Z'_k)_{k=1:n}$ instead of $(Z_k)_{k=1:n}$. The processes $(X_k)_{k=0:n}$ and $(X'_k)_{k=0:n}$ are $({\cal F}^Z_k=\sigma(X_0,(Z_\ell)_{\ell=1:k}))_{k=0:n}$  and $({\cal F}^{Z'}_k=\sigma(X'_0,(Z'_\ell)_{\ell=1:k}))_{k=0:n}$-Markov chains respectively. It is clear by backward induction that
\[
\Phi_k(X_{0:k}) = \E\big( \Phi_n(X_{0:n})\,|\, {\cal F}^Z_k) \quad\mbox{ and}\quad \Phi'_k(X_{0:k}) = \E\big( \Phi'(X'_{0:n})\,|\, {\cal F}^{Z'}_k), \; k=0, \ldots,n.
\]
Let us suppose that {\bf the assumption is satisfied by $(Z_{k+1},\vartheta_k,\vartheta'_k)_{k=0:n-1}$} : for each $k=0, \ldots,n-1$,~\eqref{eq:convmat} holds and $Z_{k+1}$ has a radial distribution. Since the law of $Z_{k+1}$ is radial, \begin{equation}
\forall (x_{0:k},u,O)\!\in (\R^d)^{k+1}\times {\mathbb M}_{d,q}\times {\mathbb M}_{q,q}\mbox{ with $O$ orthogonal }\Psi_k(x_{0:k}, u)=\Psi_k(x_{0:k}, uO)\label{eq:invgkortho}.
\end{equation}
We first check  by backward induction that the functionals $\Phi_k$ are convex. The function $\Phi_n$   is convex  by assumption. If $\Phi_{k+1}$ is convex, by convexity of $\R^d\ni w\mapsto \Phi_{k+1}(x_{0:k},x_k+w)$ and Lemma~\ref{pro:Cvxdq}$(iii)$,
\begin{equation}
   \forall x_{0:k}\in (\R^d)^{k+1},\;\forall u,v\in {\mathbb M}_{d,q}\mbox{ s.t. }uu^*\le vv^*,\;\Psi_k(x_{0:k}, u)\le \Psi_k(x_{0:k},v).\label{eqcroisG}
\end{equation}
With~\eqref{eq:convmat} then the convexity of $\Psi_{k}$ consequence of the one of $\Phi_{k+1}$ and last~\eqref{eq:invgkortho}, we deduce that for $x_{0:k},y_{0:k}\in(\R^d)^{k+1}$ and $\alpha\in[0,1]$, 
\begin{align*}
\Phi_k\big(\alpha x_{0:k}+(1-\alpha)y_{0:k}\big) &=  \Psi_k\big(\alpha x_{0:k}+(1-\alpha)y_{0:k},\vartheta_k(\alpha x_{k}+(1-\alpha)y_{k})\big)\\
                                                 &\le  \Psi_k\big(\alpha x_{0:k}+(1-\alpha)y_{0:k},\alpha \vartheta_k(x_{k})+(1-\alpha)\vartheta_k(y_{k})O_{k,x_k,y_k,\alpha}\big)\\
                                                 &\le \alpha \,\Psi_k\big(x_{0:k},\vartheta_k(x_{k})\big)+(1-\alpha) \Psi_k\big(y_{0:k},\vartheta_k(y_{k})O_{k,x_k,y_k,\alpha}\big).\\
						&= \alpha \,\Psi_k\big(x_{0:k},\vartheta_k(x_{k})\big)+(1-\alpha) \Psi_k\big(y_{0:k},\vartheta_k(y_{k})\big)\\  
						&=\alpha\,\Phi_k\big(x_{0:k}\big)+(1-\alpha)\Phi_k\big(y_{0:k}\big).
\end{align*}


\smallskip
As a second step, let us prove that $\Phi'_k\le \Phi_k$, $k=0, \ldots,n$, still by  backward induction. This is true for $k=n$ since $\Phi_n= \Phi'_n$. Assume $\Phi'_{k+1}\le \Phi_{k+1}$. Then, 
\[
\Psi'_k(x_{0:k}, u) \le \E\, \Phi_{k+1}\big(x_{0:k}, x_k+uZ'_{k+1}\big).
\] 
Now, for every $(x_{0,k},u)\!\in (\R^d)^{k+1}\times{\mathbb M}_{d,q}$, the function $z\mapsto  \Phi_{k+1}\big(x_{0:k}, x_k+u\, z\big)$ is convex as the composition of a convex function with an affine function. The assumption $Z'_{k+1}\le_{cvx}Z_{k+1}$ implies that 
\[
\Psi'_k(x_{0:k},u)\le \E\, \Phi_{k+1}\big(x_{0:k}, x_k+uZ_{k+1}\big)= \Psi_k(x_{0:k},u)
\] 
which in turn ensures, once composed with $\vartheta'_k$,  that $\Phi'_k(x_{0:k})\le \Psi_k(x_{0:k},\vartheta'_k(x_k))$. With the condition $\vartheta'_k(\vartheta_k')^*\le \vartheta_k\vartheta_k^*$ and~\eqref{eqcroisG} we deduce that $\Phi'_k\le \Phi_k$.  Since this inequality holds for every $k$, one has in particular that $\Phi'_0\le \Phi_0$ so that
\[
\E\, \Phi'_n\big(X'_{0:n}\big)=\E\, \Phi_0'(X'_0)\le \E\, \Phi_0(X'_0)\le \E\, \Phi_0(X_0)=\E\, \Phi_n\big(X_{0:n}\big),
\]
where we used in the last inequality the assumption $X'_0\le_{cvx} X_0$ and the convexity of $\Phi_0$.

When {\bf the assumption is satisfied by $(Z'_{k+1},\vartheta'_k,\vartheta_k)_{k=0:n-1}$} i.e. for each $k=0, \ldots,n-1$,~\eqref{eq:convmat} holds with $\vartheta_k'$ replacing $\vartheta_k$ and $Z'_{k+1}$ has a radial distribution, then we check as above that, for each $k=0, \ldots,n-1$, $\Phi'_k$ is convex and that
\begin{equation}
   \forall x_{0:k}\in (\R^d)^{k+1},\;\forall u,v\in {\mathbb M}_{d,q}\mbox{ s.t. }uu^*\ge vv^*,\;\Psi'_k(x_{0:k}, u)\ge \Psi'_k(x_{0:k},v).\label{eq:croisg'k}
\end{equation}
To deduce  by  backward induction that $\Phi_k\ge \Phi'_k$, $k=0, \ldots,n$, we assume $\Phi_{k+1}\ge \Phi'_{k+1}$. Then, 
\[
\Psi_k(x_{0:k}, u) \ge \E\, \Phi'_{k+1}\big(x_{0:k}, x_k+uZ_{k+1}\big)\ge \E\, \Phi'_{k+1}\big(x_{0:k}, x_k+uZ'_{k+1}\big)=\Psi'_k(x_{0:k}, u),
\]
where we used the convexity of $\Phi'_{k+1}$ and $Z'_{k+1}\le_{cvx}Z_{k+1}$ for the second inequality. By composing with $\vartheta_k$ then using~\eqref{eq:croisg'k} with $u=\vartheta_k(x_k)$ and $v=\vartheta'_k(x_k)$ thanks to the condition $\vartheta_k\vartheta_k^*\ge \vartheta'_k(\vartheta_k')^*$, we deduce that $$\Phi_k(x_{0:k})\ge \Psi'_k(x_{0:k},\vartheta_k(x_k))\ge \Psi'_k(x_{0:k},\vartheta'_k(x_k))=\Phi'_k(x_{0:k}).$$ One has in particular that $\Phi_0\ge \Phi'_0$ so that
\[
\E\, \Phi_n\big(X_{0:n}\big) = \E\, \Phi_0(X_0) \ge \E\, \Phi'_0(X_0)\ge \E\, \Phi_0'(X'_0) = \E\, \Phi'_n\big(X'_{0:n}\big) 
\]
where we used in the last inequality the assumption $X'_0\le_{cvx} X_0$ and the convexity of $\Phi_0'$.
%
%
~$\hfill\cqfd$

\bigskip

%

 \bigskip
\noindent {\bf A generalization of Proposition~\ref{Proconvgen}.}  When $q=d$ and in particular in the one-dimensional case $q=d=1$, we can relax the radial assumption on the $Z_k$ owing to  properties $(i)$ and $(ii)$ from Lemma~\ref{pro:Cvxdq} rather than on claim $(iii)$. This should be compensated by strenghtening the assumptions on  diffusion coefficients $\vartheta_k$. This leads to the following three assumptions:

\smallskip Either,  for each $k=0, \ldots,n-1$, one of the following three conditions holds  ({\bf  we say that the assumption is satisfied by $(Z_{k+1},\vartheta_k,\vartheta'_k)_{k=0:n-1}$}):
\begin{itemize}
  \item $Z_{k+1}$ satisfies the vanishing conditional expectation assumption, $\vartheta_k$ and $\vartheta'_k$ are diagonal  both with non-negative entries, the ones of $\vartheta_k$ being moreover convex  and $\vartheta_k'(\vartheta_k')^*\le \vartheta_k\vartheta_k^*$ (i.e. $(\vartheta'_k)_{ii}\le (\vartheta_k)_{ii}$, $i=1:d$),
  \item  $Z_{k+1}$ satisfies the symmetric conditional  distribution assumption, $\vartheta_k$ and $\vartheta'_k$ are diagonal with the entries of $\vartheta_k$ convex and $\vartheta_k'(\vartheta_k')^*\le \vartheta_k\vartheta_k^*$ (i.e. $|(\vartheta'_k)_{ii}|\le |(\vartheta_k)_{ii}|$, $i=1:d$),  
  \item  $Z_{k+1}$ has a radial distribution,  $\vartheta_k$ is convex in the the matrix-convexity sense~\eqref{eq:convmat}  and $\vartheta_k'(\vartheta_k')^*\le \vartheta_k\vartheta_k^*$,
  \end{itemize}
 
  \smallskip
 \noindent or  ({\bf assumption satisfied by $(Z'_{k+1},\vartheta'_k,\vartheta_k)_{k=0:n-1}$})  for each $k=0, \ldots,n-1$ in the following sense: $\vartheta_k'(\vartheta_k')^*\le \vartheta_k\vartheta_k^*$  and  $(Z_k,\vartheta_k,\vartheta'_k)$ is replaced by  $(Z'_k,\vartheta'_k,\vartheta_k)$ in  other assertions.
  

\begin{Pro}[Convex order:  $q=d$]\label{propconvdiag} Let $(Z_k)_{k=1:n}$ and $(Z'_k)_{k=1:n}$ be two sequences of independent integrable and centered $\R^d$-valued random vectors. Let $(\vartheta_k)_{k=0:n-1}$ and $(\vartheta'_k)_{k=0:n-1}$ be two sequences of ${\mathbb M}_{d,d}$-valued functions with linear growth defined on $\R^d$. Let $X_0$ and $X'_0$ be integrable $\R^d$-valued r.v. independent of $(Z_k)_{k=1:n}$ and $(Z'_k)_{k=1:n}$ respectively.
  Denote by $(X_k)_{k=0:n}$ and $(X'_k)_{k=0:n}$ the two {\em ARCH} models respectively defined by~\eqref{eq:ARCHdef} and by $X'_{k+1}=X'_{k}+\vartheta'_k(X'_{k})Z'_{k+1}$ for $k=0, \ldots,n-1$.

  Under the assumption just before the proposition and if $X'_0\le_{cvx} X_0$ and $Z'_k\le_{cvx} Z_k$ for every $k=1:n$, then 
$$
(X'_k)_{k=0:n} \le_{cvx} (X_k)_{k=0:n}.
$$
\end{Pro}
\noindent {\bf Proof.}
The result is a special case of Proposition~\ref{Proconvgen} when $(Z_k)_{k=1:n}$ satisfies the radial distribution assumption and $\vartheta_k'(\vartheta_k')^*\le \vartheta_k\vartheta_k^*$. We are simply going to explain how to adapt the backward induction steps in the proof of this proposition when the assumption before the proposition is satisfied by $(Z_{k+1},\vartheta_k,\vartheta'_k)_{k=0:n-1}$ and $Z_{k+1}$ satisfies either the vanishing conditional expectations assumption or the symmetric conditional distribution assumption and the matrices  $\vartheta_k$ and $\vartheta_k'$ are diagonal.
Let $\Phi_n: (\R^d)^{n+1}\to \R$ be a convex function with linear growth and $\Phi'_n=\Phi_n$. We define by backward induction the sequence $(\Psi_k,\Phi_k,\Psi'_k,\Phi'_k)_{k=0:n-1}$ by 
  \begin{align*}
    &\Psi_k(x_{0:k},u)=\E\,\Phi_{k+1}(x_{0:k},x_k+{\rm Diag}(u)Z_{k+1})\mbox{ and  }\Phi_{k}(x_{0:k})=\Psi_k(x_{0:k},(\vartheta_k(x_k)\mbox{\bf 1}))\\
    &\Psi'_k(x_{0:k},u)=\E\,\Phi'_{k+1}(x_{0:k},x_k+{\rm Diag}(u)Z'_{k+1})\mbox{ and  }\Phi'_{k}(x_{0:k})=\Psi'_k(x_{0:k},(\vartheta'_k(x_k)\mbox{\bf 1}))
  \end{align*}
  where $x_{0:k}\in(\R^d)^{k+1}$, $u\in\R^d$, ${\rm Diag}(u)\in{\mathbb M}_{d,d}$ denotes the diagonal matrix with diagonal coefficients ${\rm Diag}(u)_{ii}=u_i$, $i=1:d$, and $\mbox{\bf 1}$ is the vector in $\R^d$ with all coefficients equal to $1$. Note that when $\vartheta_k$ is diagonal ${\rm Diag}(\vartheta_k(x)\mbox{\bf 1})=\vartheta_k(x)$ for all $x\in\R^d$.
If $\Phi_{k+1}$ is convex and $Z_{k+1}$ satisfies the vanishing conditional expectations assumption (or the stronger symmetric conditional distributions assumption), then by the convexity of $\R^d\ni w\mapsto \Phi_{k+1}(x_{0:k},x_k+w)$ and Lemma~\ref{pro:Cvxdq}$(i)$,
\begin{equation}
   \forall x_{0:k}\in (\R^d)^{k+1},\;\forall u,v\in \R_+^d\mbox{ s.t. }u_i\le v_i\mbox{ for }i=1:d,\;\Psi_k(x_{0:k}, u)\le \Psi_k(x_{0:k},v).\label{eq:croisGbis}
\end{equation}
For $u\in\R^d$, let us denote by ${\rm abs}(u)$ the vector in $\R^d$ defined by ${\rm abs}(u)_i=|u_i|$, $i=1:d$. Assume moreover either that $\vartheta_k$ is diagonal with nonnegative and convex diagonal coefficients or $Z_{k+1}$ satisfies the symmetric conditional distributions assumption, $\vartheta_k$ is diagonal and the absolute values of its diagonal coefficients are convex functions. Then \begin{equation}\label{eq:egalpsiabs}
   \Psi_k(x_{0:k},{\rm abs}(\vartheta_k(x_k)\mbox{\bf 1}))=\Psi_k(x_{0:k},\vartheta_k(x_k)\mbox{\bf 1})=\Phi_k(x_{0:k})
\end{equation} with the coefficients of ${\rm abs}(\vartheta_k(x_k)\mbox{\bf 1})$ nonnegative and convex. With~\eqref{eq:croisGbis} then the convexity of $\Psi_{k}$ consequence of the one of $\Phi_{k+1}$, we deduce that for $x_{0:k},y_{0:k}\in(\R^d)^{k+1}$ and $\alpha\in[0,1]$, 
\begin{align*}
\Phi_k\big(\alpha x_{0:k}+(1-\alpha)y_{0:k}\big) &=  \Psi_k\big(\alpha x_{0:k}+(1-\alpha)y_{0:k},{\rm abs}(\vartheta_k(\alpha x_{k}+(1-\alpha)y_{k})\mbox{\bf 1})\big)\\
                                                 &\le  \Psi_k\big(\alpha x_{0:k}+(1-\alpha)y_{0:k},\alpha\, {\rm abs}(\vartheta_k(x_{k})\mbox{\bf 1})+(1-\alpha){\rm abs}(\vartheta_k(y_{k})\mbox{\bf 1})\big)\\
                                                 &\le \alpha \Psi_k\big(x_{0:k},{\rm abs}(\vartheta_k(x_{k})\mbox{\bf 1})\big)+(1-\alpha) \Psi_k\big(y_{0:k},{\rm abs}(\vartheta_k(x_{k})\mbox{\bf 1})\big)\\
  &=\alpha\Phi_k\big(x_{0:k}\big)+(1-\alpha)\Phi_k\big(y_{0:k}\big).
\end{align*}
If $\Phi'_{k+1}\le \Phi_{k+1}$, then, 
\[
\Psi'_k(x_{0:k}, u) \le \E\, \Phi_{k+1}\big(x_{0:k}, x_k+{\rm Diag}(u)Z'_{k+1}\big).
\] 
Now, for every $(x_{0,k},u)\!\in (\R^d)^{k+1}\times\R^d$, the function $z\mapsto  \Phi_{k+1}\big(x_{0:k}, x_k+{\rm Diag}(u)\, z\big)$ is convex as the composition of a convex function with an affine function. The assumption $Z'_{k+1}\le_{cvx}Z_{k+1}$ implies that 
\[
\Psi'_k(x_{0:k},u)\le \E\, \Phi_{k+1}\big(x_{0:k}, x_k+{\rm Diag}(u)Z_{k+1}\big)= \Psi_k(x_{0:k},u)
\] 
which in turn ensures, once composed with $\vartheta'_k$,  that $$\Phi'_k(x_{0:k})\le \Psi_k(x_{0:k},\vartheta'_k(x_k)\mbox{\bf 1})=\Psi_k(x_{0:k},{\rm abs}(\vartheta'_k(x_k)\mbox{\bf 1})).$$ Since the absolute values of the diagonal coefficients of $\vartheta'_k$ are not greater than the ones of $\vartheta_k$, we deduce with~\eqref{eq:croisGbis} and \eqref{eq:egalpsiabs} that $\Phi'_k\le \Phi_k$.
~$\hfill\cqfd$
\bigskip

In the scalar $q=d=1$ case, we deduce the following result.
\begin{Pro}[Scalar setting: $d=q=1$]\label{propconvdomd1} Let $(X_k)_{k=0:n}$be a scalar {\em ARCH} model defined by~\eqref{eq:ARCHdef} where the white noise $(Z_k)_{k=1:n}$ is  scalar but (possibly) not bounded and $g_0(X_0)\le_{cvx} X_0$ with $X_0$ integrable. Assume that the functions $|\vartheta_k|$, $k=0,\ldots,n-1$ are  convex with linear growth and that for each $k=0, \ldots,n-1$, either $\vartheta_k$ is nonnegative or $-Z_{k+1}$ has the same distribution as $Z_{k+1}$.  

\smallskip
\noindent $(a)$ {\em Truncation}. Let  $(\breve Z_k=Z^{A_k}_k)_{k=1:n}$ with $(A_k)_{k=1:n}$
an $n$-tuple of Borel sets satisfying 
$ \E\, Z_k\mbox{\bf 1}_{\{Z_k \in A_k\}}=0$, $k=1,\ldots,n$ and  $(\breve X^A_k)_{k=0:n}$  be  the induced approximating {\em ARCH} process defined by~\eqref{eq:ARCHTruncdef}. Then
\[
(\breve X^A_k)_{k=0:n} \le_{cvx} (X_k)_{k=0:n}.
\]  

\noindent $(b)$ {\em Voronoi  Quantization}. Let $(\breve Z_{k}= \widehat Z^{vor}_k)_{k=1:n}$ be a stationary (Voronoi) quantized approximation of the white noise $Z_{1:n}$ and $\breve X_{0:n}$ the induced approximating {\em ARCH} process defined by~\eqref{eq:ARCHTruncdef}. Then 
\[
(\breve X_k)_{k=0:n} \le_{cvx} (X_k)_{k=0:n}.
\]
\end{Pro}
When $g_0$ is  nearest neighbour projection on a stationary Voronoi (primal) quantization grid for $X_0$, the hypothesis $g_0(X_0)\le_{cvx} X_0$ is satisfied.

\bigskip
\noindent {\bf Proof.} $(a)$ Follows from the combination of Lemma~\ref{lem:truncZ} and Proposition~\ref{propconvdiag}.

\smallskip
\noindent $(b)$ Follows from the stationarity property which implies $ \widehat Z^{vor}_k = \E\big(Z_k\,|\,  \widehat Z^{vor}_k\big) \le_{cvx} Z_k$, $k=1:n$ and Proposition~\ref{propconvdiag}.~$\hfill\cqfd$

\subsection{Convex domination of $(\breve X_k)_{k=0:n}$ by $(\widehat X_k)_{k=0:n}$}\label{secpathconv}
\begin{Thm}\label{thmdomgausd} Assume that the ${\mathbb M}_{d,q}$-valued functions $\vartheta_k$, $k=0,\ldots,n-1$, are convex in the sense of~\eqref{eq:convmat} with linear growth and the random vectors $\breve Z_k$, $k=1:n$ have a radial distribution or that $q=d$ and the assumption before Proposition \ref{propconvdiag} is satisfied by $(\breve Z_{k+1},\vartheta_k,\vartheta_k)_{k=0:n-1}$. If $\breve X_0\le_{cvx}\widehat X_0$, then the processes $(\breve X_k)_{k=0:n}$ and $(\widehat X_k)_{k=0:n}$ respectively defined by \eqref{eq:ARCHTruncdef} and \eqref{eq:XZtildebreve} satisfy $$(\breve X_k)_{k=0:n}\le_{cvx}(\widehat X_k)_{k=0:n}.$$
\end{Thm}
\begin{Remark} When $Z_k$ has a radial distribution, then this property is inherited by the truncated noise $\breve Z_k=Z_k{\mathbf 1}_{Z_k\in A_k}$ as soon as the set $A_k$ is invariant by orthogonal transformations and in particular if it is an Euclidean ball centered at the origin. By contrast, the radial distribution is not in general preserved through optimal quadratic primal quantization. If $Z_k$ has independent coordinates, then $\breve Z_k$ obtained by product quantization or truncation outside product sets also has independent coordinates and therefore satisfies the vanishing conditional expectation condition. 
\end{Remark}

\noindent {\bf Proof.}
Let $\Phi_n: (\R^d)^{n+1}\to \R$ convex with linear growth. We adapt the proof of Proposition \ref{Proconvgen} or, when $q=d$, of Proposition \ref{propconvdiag} by defining by backward induction starting from $\Phi_n=\Phi'_n$ :
\begin{align*}
  &\Psi_k(x_{0:k},u)=\E\,\Phi_{k+1}(x_{0:k},{\rm Proj}^{del}_{\Gamma_{k+1}}(x_k+u\breve Z_{k+1}, U_{k+1})),\;\Phi_k(x_{0:k})=\Psi_k(x_{0:k},\vartheta_k(x_k))\\
  &\Psi'_k(x_{0:k},u)=\E\,\Phi'_{k+1}(x_{0:k},x_k+u\breve Z_{k+1}),\;\Phi'_k(x_{0:k})=\Psi'_k(x_{0:k},\vartheta_k(x_k)),\;k=0, \ldots,n-1,
\end{align*}
Since, except in the scalar case $d=1$, the convex property is not necessarily preserved through the dual quantization step, the convexity of the functions $\Phi_k$ is not clear. Nevertheless, arguing like in the proof of Proposition \ref{Proconvgen} or Proposition \ref{propconvdiag}, we check that the functions $\Phi'_k,\;k=0:n$ are convex. Let us deduce  that $\Phi_k\ge \Phi'_k$ be backward induction on $k=0:n$. Equality holds for $k=n$ and supposing that $\Phi_{k+1}\ge \Phi'_{k+1}$, we have, using Jensen's inequality combined with the stationarity property \eqref{eq:statprop} for the second inequality,
\begin{align*}
  \Psi_k(x_{0:k},u)&\ge \E\Big[\E\big(\Phi'_{k+1}(x_{0:k},{\rm Proj}^{del}_{\Gamma_{k+1}}(x_k+u\breve Z_{k+1}, U_{k+1}))|\breve Z_{k+1}\big)\Big]\\
  &\ge \E\,\Phi'_{k+1}(x_{0:k},x_k+u\breve Z_{k+1})=\Psi'_k(x_{0:k},u)
\end{align*}
so that for the choice $u=\vartheta_k(x_k)$, we conclude that $\Phi_k\ge \Phi'_k$. Hence, using the convexity of $\Phi_0$ and $\breve X_0\le_{cvx}\widehat X_0$ for the second inequality, we have
$$
\hskip 2.5cm \E\,\Phi_n(\widehat X_{0:n})=\E\,\Phi_0(\widehat X_0)\ge \E\,\Phi'_0(\widehat X_0)
\ge \E\,\Phi'_0(\breve X_0)=\E\,\Phi_n(\breve X_{0:n}).\hskip2.5cm\cqfd
$$

\section{Application to the Euler scheme of a martingale Brownian diffusion} \label{sec:euler}
The   {\em Euler scheme} with step $h = \tfrac Tn$, $n\ge 1$, and Brownian increments of a martingale Brownian diffusion $dX_t=\vartheta(t,X_t)dW_t$ is an {\em ARCH} model corresponding  to the choices: $\vartheta_k(x)= \sqrt{\frac Tn}\vartheta(t_k, x)$ with $t_k= \tfrac{kT}{n}$ and $(Z_k=\sqrt{\tfrac nT }(W_{t_k}-W_{t_{k-1}}))_{k=0:n}$ a Gaussian  ${\cal N}(0;I_q)$-distributed white noise ($W$ is a standard $q$-dimensional Brownian motion). It reads
\begin{equation}\label{eq:Euler}
\bar X_{k+1} = \bar X_k +\vartheta(t_k, \bar X_k) \sqrt{\tfrac Tn } Z_{k+1}, \; k=0, \ldots,n-1.
\end{equation}


We can consider three different ways to modify the white noise:

\smallskip
-- Replacing the normalized Brownian increments  by a $q$-dimensional Rademacher white noise, one drawback being that we loose all pathwise comparison with the original diffusion combined with an explosion  of the state space (of size $2^q$) in comparison with primal quantization  when the dimension $q$ of the white noise is large. 

\smallskip
-- Truncating (see below) the Gaussian white noise   may be convenient in higher dimension, the drawback being this time that the distribution of the approximation $(\widehat X_k)_{k=0:n}$ will require some Monte Carlo simulation of $(\breve X_k)_{k=0:n})$  to evaluate the transitions matrices of the chain $(\widehat X_k)_{k=0:n}$.

\smallskip
-- Quantizing the white noise by setting  $\breve{Z}_k = \widehat Z_k^{vor}= {\rm Proj}_{\Gamma^Z_k}(Z_k)$ $k=1:n$ (see~\eqref{eq:VorQuantZ} in Section~\ref{subsec:Qscheme}).  This is possible when this quantization can be performed optimally by computational means like in one  ($d=q=1$, see~\cite{Pag2015}) or   medium dimensions, in fact  as long as optimal quadratic Voronoi quantization can be computed by deterministic algorithms (see Section \ref{App:A1} for dimension one or~\cite{Montes19} for dimension $2$). For Gaussian white noise, pre-computed $L^2$-optimal grids of the ${\cal N}(0;I_q)$-distributions are available for download on the website {\tt www.quantize.maths-fi.com} up to dimension $q=10$ and size $N=1\,500$. In such settings,  both the optimization of the grids (Voronoi and dual) and the computations of the transitions matrices of the chain $(\widehat X_k)_{k=0:n}$ is extremely fast  as emphasized by recent works on Markovian recursive quantization (see~\cite{PaSa1, PaSa2} and the references therein) which is close in spirit to our primal/dual recursive quantization procedure.

\subsection{Truncation or quantization}
\paragraph{$\blacktriangleright$  Truncation of the Euler scheme with Gaussian increments\\}Assume that the diffusion function $\vartheta(t,x)$ is Lipschitz continuous in $x$ with constant $[\vartheta]_{\rm Lip}$ uniformly in $t\!\in[0,T]$ and set $c(\vartheta)= \sup_{(t,x)\in [0,T]\times\R^d} \frac{\vertiii{\vartheta(t,x)}^2}{1+|x|^2}$, $c_{\rm Fr}(\vartheta)= \sup_{(t,x)\in [0,T]\times\R^d} \frac{\|\vartheta(t,x)\|^2_{\rm Fr}}{1+|x|^2}$ and $C(\vartheta)=(q[\vartheta]^2_{\rm Lip})\vee c_{\rm Fr}(\vartheta)$. Then $\vartheta_k(x)= \sqrt{\frac{T}{n}} \,\vartheta(\frac{kT}{n},x)$, $k=0, \ldots,n-1$ and
\begin{align}
   &\max_{0\le k\le n-1}[\vartheta_k]^2_{\rm Lip} \le \frac Tn[\vartheta]^2_{\rm Lip},\;\max_{0\le k\le n-1}c(\vartheta_k)\le \frac Tn c(\vartheta),\;\notag\\&\max_{0\le k\le n-1}c_{\rm Fr}(\vartheta_k)\le {\frac Tn}c_{\rm Fr}(\vartheta),\;\max_{0\le k\le n-1}C(\vartheta_k) \le {\frac Tn} C(\vartheta).\label{eq:scalconst}
\end{align}
With Proposition \ref{prop:3.1}, we deduce that
\begin{align*}
 \|X_{k}-\breve X_k\|_2^2& \le  \| X_0-\breve X_0\|_2^2 \Big(1+ q\frac Tn[\vartheta]_{\rm Lip}^2\Big)^{\!k}+\big(1+\|X_0\|^2_2\big)\Big(1+\frac TnC(\vartheta)\Big)^{\!k} \frac Tn c(\vartheta) \sum_{\ell=1}^{k} \big\|Z_\ell\mbox{\bf 1}_{\{Z_\ell\in A^c_{\ell}\}}\big\|_2^2, \\
 & \le   \| X_0-\breve X_0\|_2^2 e^{q[\vartheta]_{\rm Lip}^2\frac{kT}{n}}+\big(1+\|X_0\|^2_2\big)e^{C(\vartheta) \frac{kT}{n}}c(\vartheta)\frac{kT}{n}\big\|Z_1\mbox{\bf 1}_{\{|Z_1|\ge a\}}\big\|_2^2, \quad k=0, \ldots,n,
\end{align*}
so that, by Doob's Inequality, 
\[
\Big\|  \max_{k=0:n}|X_{k}-\breve X_k| \Big\|_2^2\le 4\| X_0-\breve X_0\|_2^2 e^{q[\vartheta]_{\rm Lip}^2T}+4 T  \big(1+\|X_0\|^2_2\big)e^{C(\vartheta) T}c(\vartheta)\big\|Z_1\mbox{\bf 1}_{\{|Z_1|\ge a\}}\big\|_2^2.
\]
\noindent  {\em Choice of $a= a(n)$.} -- If $q=1$, the tail expectation can be estimated by a straightforward integration by parts which: for every $a>0$
\[
\E \, |Z_1|^2\mbox{\bf 1}_{\{|Z_1|\ge a\}} \le \sqrt{\frac{2}{\pi}}\Big(a+\frac 1a\Big) e^{-\frac{a^2}{2}},\quad a > 0.
\]
If we set $a = a_n \ge   \sqrt{c\, \log n}$ for some  $c>0$, then  
\[
\E \, |Z_1|^2\mbox{\bf 1}_{\{|Z_1|\ge a\}} = O\left(\frac{\sqrt{\log n}}{n^{\frac c2 }}\right) \to 0.
\]
-- If $q\ge 2$, a simple, though sub-optimal, approach is the following: we start form the obvious
\begin{align*}
\E \, |Z_1|^2\mbox{\bf 1}_{\{|Z_1|\ge a\}}& \le e^{-\lambda a^2} \E \,|Z_1|^2e^{\lambda |Z_1|^2} =  e^{-\lambda a^2}q \,\E\, \zeta^2e^{\lambda \zeta^2}\times  \Big(\E\,  e^{\lambda \zeta^2}\Big)^{d-1}\\
&= q\,\frac{ e^{-\lambda a^2} }{(1-2\lambda)^{ \frac{3}{2}+\frac{d-1}{2} } }\quad \mbox{ where }\quad \zeta \sim {\cal N}(0;1).
\end{align*}

As soon as $a> \sqrt{d+2}$,  the function $\lambda\mapsto -\lambda a^2 -\frac{d+2}{2} \log(1-2\lambda)$ is minimum at $\lambda(a) = \frac 12 \Big(1-\frac{d+2}{a^2}\Big)\!\in (0, \frac 12)$. Hence
\[
\E \, |Z_1|^2\mbox{\bf 1}_{\{|Z_1|\ge a\}} \le q\,e^{-\frac{a^2}{2}}\left( \frac{e a^2}{d+2}\right)^{1+\frac d2}
\]
and, if  $a = a_n \ge   \sqrt{c\, \log n}$ for some  $c>0$, then  
\[
\E \, |Z_1|^2\mbox{\bf 1}_{\{|Z_1|\ge a\}} = O\left(\frac{ (\log n)^{1+\frac d2}}{n^{\frac c2 }}\right) \to 0.
\]
\paragraph{$\blacktriangleright$ Voronoi/primal quantization of the increments of the Euler scheme\\} As ${\cal N}(0;I_q)$ has $2+\eta$-moment for any $\eta>0$, it follows from Zador's Theorem (see Theorem~\ref{thm:zadoretpierce1}) that, if $\breve Z_k$ are either  quadratic optimal quantizations of $Z_k$ at level $N^Z_k$ or (like  in the above remark) a product quantization of optimal quantizations of the marginal, in both cases 
\[
\big\| Z_k  - \breve Z_k\|_2 = e_{2,N_Z}\big({\cal N}(0;I_q))\le C_{q,\eta} \sigma_{2+\eta}({\cal N}(0;I_q))(N^Z_k)^{-1/d},
\] 
where $\sigma_{2+\eta}({\cal N}(0;I_q))=  \Big(\frac{2^{\eta/2}}{\pi^{q/2}} S_{q-1}\Gamma\big(\frac{\eta+q}{2}+1\big)\Big)^{\frac{1}{2+\eta}}$ with $\Gamma(.)$ denoting the Euler $\Gamma$ function and $ S_{q-1}$ the area of the unit  sphere of dimension $q-1$.

\subsection{Approximation rate} 

If we assume that the diffusion function $\vartheta(t,x)$ is Lipschitz in $x$ uniformly in $t\!\in [0,T]$, then, under the assumptions of the above claim~$(b)$ of Theorem~\ref{thm:cvgce}, we obtain using \eqref{eq:scalconst} that for every $k=0, \ldots,n$,
\begin{align*}
 \big \|\widehat X_{k}-\bar X_{k}&\big\|_2 \le \left( (\widetilde C^{vor}_{2,d,\eta})^2e^{qt_k[\vartheta]^2_{\rm Lip} }  \frac{\sigma^2_{2+\eta}(X_0)}{N_0^{2/d}}\right.
\\
 & \quad \left.+\sum_{\ell=1}^ke^{q(t_k-t_\ell)[\vartheta]^2_{\rm Lip} }   \left[\frac Tn  \|\vartheta (t_{\ell-1}, \bar X_{\ell-1})\|_2^2(\widetilde C^{vor}_{2,q,\eta})^2\frac{\sigma_{2+\eta}^2(Z_\ell)}{(N^Z_\ell)^{2/q}} + (\widetilde C^{del}_{2,d,\eta})^2 \frac{\sigma_{2+\eta}^2(\widetilde X_{\ell})}{(N_\ell)^{2/d}}\right] \right)^{1/2}.
\end{align*}
Note that,  according to~\eqref{eq:Bound0},
\[
 \|\vartheta (t_{k}, \bar X_{k})\|_2^2\le c(\vartheta)e^{c_{\rm Fr}(\vartheta) t_k}(1+\|X_0\|_2^2). 
\]

\subsection{Domination}

If the diffusion coefficient $\vartheta: [0,T]\times \R^d \to \mathbb{M}_{d,q}(\R)(t,x)$ satisfies $\vartheta(t,\cdot)$ is convex in the sense of~\eqref{eq:convmat} and has linear growth in $x$ uniformly if $t\!\in [0,T]$, then the {\em ARCH} $(\breve X_k)_{k=0:n}$ with either truncated  (by ball $B(0;a)$) or optimally Voronoi quantized noise (but not with Rademacher white noise) is dominated for the convex order by the original Euler scheme with Gaussian white noise according to Theorem \ref{propdomgausd}. This theorem does not apply when $(\breve X_k)_{k=0:n}$ is driven by a Rademacher white noise. Indeed, the $q$-dimensional Rademacher distribution is not comparable to the standard Gaussian ${\cal N}(0;I_q)$ distribution for the convex order since for $i\in\{1,\cdots,q\}$, the integral of the convex function $\R^q\ni (x_1,\cdots,x_q)\mapsto |x_i|$ with respect to the former (resp. latter) is equal to $1$ (resp. $\sqrt{2/\pi}<1$) while the integral of $(x_i)^4$ with respect to the former (resp. latter) is equal to $1$ (resp. $3>1$). Last, by Theorem \ref{thmdomgausd}, the {\em ARCH} $(\breve X_k)_{k=0:n}$ with truncated by ball $B(0;a)$ noise is dominated by the scheme $(\widehat X_k)_{k=0:n}$where a dual quantization step follows the truncated {\em ARCH} dynamics at each time, since the truncated noise is radially distributed.

\providecommand{\norm}[1]{\lVert#1\rVert}
\providecommand{\pnorm}[1]{\norm{#1}_p}
\providecommand{\abs}[1]{\lvert#1\rvert}
\providecommand{\enorm}[1]{\abs{#1}}
\providecommand{\ind}[1]{\mathbbm{1}_{#1}}
\providecommand{\card}[1]{\cardop{#1}}
\newcommand{\Qn}{Q_{\text{min}}}
\newcommand{\Qx}{Q_{\text{max}}}
\newcommand{\Qkn}{Q^k_{\text{min}}}
\newcommand{\Qkx}{Q^k_{\text{max}}}
\newcommand{\qXk}{\widehat X_k}
\newcommand{\qXkk}{\widehat X_{k+1}}
\newcommand{\tp}{\pi^k_{ij}}

\appendix
\section{Background on (optimal) primal and dual vector quantization}\label{App:A}

%

In what follows $\R^d$ is supposed to be equipped with the canonical Euclidean norm. For a more general presentation dealing with any norm, see~\cite{GrLu} for Voronoi quantization and~\cite{PaWi1} for Delaunay quantization. 
\subsection{Optimal Voronoi  quantization (primal)} \label{App:A1}
Let $\Gamma= \{x_1, \ldots,x_{_N}\}\subset \R^d$ denote a finite subset of size $N$, that we will call {\em grid}. To such a grid we can associate {\em Voronoi diagrams} $(C_i(\Gamma))_{i=1:N}$ that are Borel partitions of $\R^d$ satisfying
\[
\forall\, i\!\in \{1, \ldots,N\},\quad C_i(\Gamma) \subset \big\{\xi\!\in \R^d: |\xi-x_i| \le \min_{1\le j\le N} |\xi-x_j|\big\}.
\]  
There is a one-to-one correspondence between Voronoi diagrams and Borel {\em nearest neighbour projections}, denoted ${\rm Proj}_{_\Gamma}$, defined as Borel mappings from $\R^d \to \Gamma$ such that 
\[
\forall \, \xi\!\in \R^d, \quad |\xi-{\rm Proj}_{_\Gamma}(\xi)| = {\rm dist}(\xi, \Gamma).
\]
Indeed,  if ${\rm Proj}_{_\Gamma}$ is a Borel nearest neighbour projection, then  $\big(\{{\rm Proj}_{_\Gamma}= x_i\}\big)_{i=1:N}$ is a Voronoi diagram and, conversely, for any Voronoi diagram $(C_i(\Gamma))_{i=1:N}$, 
\begin{equation}\label{eq:ProjNN}
{\rm Proj}_{_\Gamma} = \sum_{i=1}^N x_i\mbox{\bf 1}_{C_i(\Gamma)}
\end{equation}
is a Borel nearest neighbour projection. The elements $(C_i(\Gamma))$ of a Voronoi diagram are called {\em Voronoi cells}. 

We define a {\em Voronoi or primal $\Gamma$-quantization } of an $\R^d$-valued random vector $X:(\Omega, {\cal A}, \P)\to \R^d$ by 
\begin{equation}\label{eq:GammaQuantiz}
\widehat X = \widehat X^{\Gamma}:= {\rm Proj}_{_\Gamma} (X)
\end{equation}
whose distribution  is given by $\widehat \mu^{\Gamma} = \mu \circ {\rm Proj}_{_\Gamma}^{-1}$ if $X$ is $\mu$-distributed. 

If $\mu\Big(\bigcup \partial C_i(\Gamma)\Big)= 0$, then $\widehat \mu^{\Gamma}$ is unique and all $\Gamma$-quantizations are $\P$-$a.s.$ equal. The mean $L^p$-quantization error induced by $\Gamma$ is defined by 
\[
e_p(\Gamma, \mu) = e_p(\Gamma, X) = \big\| {\rm  dist}(X, \Gamma)\big\|_p= \big\| X-\widehat X^{\Gamma}\big\|_p
\]
for any Voronoi quantization of $X$ (still $\mu$-distributed). 

\smallskip
Then one defines, for $p>0$ and  a  $ N\!\in \N$,  the minimal mean $L^p$-quantization error at {\em level} $N$ by
\[
e_{p, N}(\mu)= e_{p, N}(X) = \inf_{\Gamma: |\Gamma|\le N} e_p(\Gamma, X).
\]

If $\mu$ has a finite $p$th moment, then {\em the above infimum is in fact a minimum} and any  optimal grid $\Gamma^{(N)}$ solution to the above minimization problem has a full size $N$ provided the support of $\mu$ has at least $N$ elements (see~e.g. Theorem 4.12 in~\cite{GrLu} or Theorem~5.1~in~\cite{PagSpring2018} among others). The random vector  $\widehat X^N = \widehat X^{\Gamma^{(N)}}$ is called an optimal $L^p$-quantization of $X$. Moreover, the optimal quantization $\widehat X^{N}$  is $\P$-$a.s.$ uniquely defined since one always has $\mu\Big(\bigcup \partial C_i(\Gamma^{(N)})\Big)= 0$ (see Theorem~4.2 in~\cite{GrLu}).

\smallskip Finally, in the quadratic case $p=2$, any optimal quantization grid $\Gamma^{(N)}$ at level $N$ and the associated quantization $\widehat X^N$ satisfy (see e.g.~\cite{GrLu},~\cite{Pag2015} or~\cite{PagSpring2018}, Proposition~5.1 among others)  a stationarity (or self-consistency) equation reading
\begin{equation}\label{eq:StatioVoro}
\widehat X^N = \E\,\big( X \, |\, \widehat X^N). 
\end{equation}

The key property of Optimal quantization theory is the following theorem that rules the  rate of decay of the quantization error to $0$.
\begin{Thm}[Optimal Voronoi Quantization rates]\label{thm:zadoretpierce1} 

\noindent $(a)$ \noindent {\sc Sharp rate: Zador's Theorem}, see~\cite{GrLu}: Let $X\!\in L_{\R^d}^{p+\eta}(\Omega,{\cal A}, \P)$, $p,\eta>0$,  be a random vector with distribution $\P_{_{\!X}}= \varphi.\lambda_d\stackrel{\perp}{+}\nu_{_{\!X}}$ where $\lambda_d$ denotes the Lebesgue measure and $\nu_{_X}$  denotes the singular part of the distribution. Then
\[
\lim_{N\to+\infty} N^{\frac 1d} e_{p,N}(X) = \widetilde J^{vor}_{d,p}\left(\int_{\R^d} \varphi^{\frac{d}{d+p}}d\lambda_d\right)^{\frac 1d +\frac 1p}
\]
where $\widetilde J^{vor}_{d,p}=\inf_{N\ge 1} N^{\frac 1d} d_{p,N}\big(\mathcal{U}([0,1]^d)\big)$.
%
When $d=1$, $\widetilde J^{vor}_{1,p}=\frac{1}{2(p+1)^{1/p}}$.

\medskip
\noindent $(b)$ {\sc Non-asympotic  bound (Pierce lemma)}, see~\cite{GrLu, LuPa, PagSpring2018}: Let $p, \eta >0$. For every dimension $d\ge 1$, there exists a real constant $\widetilde C^{vor}_{d,\eta,p} >0$ such that, for anyy random vector $X:(\Omega,{\cal A}, \P) \to \R^d$, 
\[
e_{p,N}(X)\le \widetilde C^{vor}_{d,\eta,p} N^{-\frac 1d} \sigma_{p+\eta}(X)
\]
where, for every $r>0$, $\sigma_r(X)= \inf_{a\in \R^d}\|X-a\|_r\le +\infty$.
\end{Thm}  
{\bf Remark.} Note that if we consider quadratic optimal {\em product}  quantizations at levels $N\ge 1$, that is solutions~--~which exist~--~to the minimization problems
\[   
e^{prod}_{2, N}(\mu)= e^{prod}_{2, N}(X) = \inf\big\{ e_2(\Gamma, X), \, \Gamma= \Gamma^1\times \cdots\times\Gamma^d, \, |\Gamma|\le N\big\}, \; N\ge 1,
\]
then, such optimal product grids are still rate optimal and satisfy a universal  non-asymptotic Pierce bound, see~e.g.~\cite{PagSpring2018}.

\medskip
Numerical computation of optimal Voronoi quantization grids of a probability distribution $\mu$ has given raise to an extensive literature starting in the 1950's originally motivated by Signal processing (see~\cite{Kieff, Trushkin}) and then developed for automatic classification (unsupervised learning) in data analysis  in close connection with $k$-means algorithms, see~\cite{McQueeen}, or more recently as a numerical method (see~\cite{Pag2015} for a survey). These methods can be divided into two classes: fixed point methods derived form the seminal Lloyd algorithm and gradient descent methods, typically the Competitive Learning Vector Quantization ($CLVQ$, see e.g.~\cite{PagSpring2018} and the references therein) with for the two families deterministic versions, quite efficient  in low dimensions ($d\le 2$ or at most $3$) and stochastic avatars for higher dimensions. 

Note that in a one-dimensional setting, if the distribution $\mu$ has a non-piecewise affine $\log$-concave density, then it is proved in~\cite{Kieff, Trushkin} that, $L^p$-optimal Voronoi  quantization grids are unique. In contrast, no such result exist in higher dimension which may have significant consequences on the efficiency of the grid optimization methods mentionned above. 

As for quadratic optimal grids of the Normal distributions in dimensions $d$ up to $10$ and sizes $N$ up to $1\, 000$ we refer to the website 

\medskip
\centerline{\tt www.quantize.maths-fi.com}

\medskip
\noindent where such grids can be downloaded (with their companion parameters like probability weights, etc).

\subsection{Optimal Delaunay (dual) quantization} \label{App:A2}

Let $X:(\Omega, {\cal A}, \P)\to \R^d$ be a random vector lying in $L^{\infty}(\P)$.  We will assume for convenience  in what follows that the support of its distribution $\mu = \P_{_X}$ spans $\R^d$ as an affine space. We can always reduce ourselves to this framework by considering the affine $A_{\mu}$ space spanned by ${\rm supp}(\mu)$ and applying  a   change of coordinates into an appropriate   orthonormal affine  basis whose $d'+1$ first points spans $A_{\mu}$ (with $d'\le d$). Optimal dual (or Delaunay) quantization relies on the best
approximation which can be achieved by a discrete random vector $\widehat X$ that
satisfies a certain stationarity assumption on the extended probability space
$(\Omega\times [0,1], {\cal A}\otimes {\cal B}([0,1]), \P\otimes \lambda_{|[0,1]})$ with $([0,1],{\cal B}([0,1]),\lambda_{|[0,1]})$ supporting a random variable uniformly distributed on $[0,1]$. 
That is why we define, for $p\!\in [1,+\infty)$:  
\begin{eqnarray*}
\forall\, N\ge d+1, \quad d_{p,N}(X) & = &  \inf_{ \widehat X}\Big\{ \big\|X -  \widehat X\big\|_p: \widehat X:(\Omega\times [0,1], {\cal A}\otimes {\cal B}([0,1]),\P\otimes \lambda_{|[0,1]})\to
\R^d,  \\
& & \qquad\qquad\qquad\qquad{\rm card}{\widehat
X(\Omega\times[0,1])} \leq N \text{ and } \E(\widehat X|X) = X \Big\}.
\end{eqnarray*}

One checks that $d_{p,N}(X)$ only depends on the distribution $\mu$ of $X$ and can subsequently be denoted $d_{p,N}(\mu)$. One shows (see~\cite{PaWi1}) that, for a given distribution $\mu$ on $(\R^d, {\cal B}or(\R^d))$,  
\begin{equation}\label{eq:dualCharac2}
d_{p,N}(\mu) = \inf\Big\{\|\Xi-\xi\|_p, \;(\Xi,\xi):(\Omega,{\cal A}, \P)\to \R^d\times\R^d, \; \Xi\sim \mu,\; \E(\xi\,|\,\Xi)=\Xi,\;  {\rm card}\big(\xi(\Omega)\big) \leq N \Big\}.
\end{equation}
Then (see~\cite{PaWi1}), one may show that such a definition is equivalent to 
\[
	d_{p,N}(X) = \inf \bigl\{ \big\|\Delta_p(X; \Gamma)\big\|_p: {\rm conv}\big({\rm supp}(\mu)\big)\subset  \Gamma \subset \R^d,
	{\rm card}(\Gamma)\leq N \bigr\}
\]
where the {\em local dual quantization functional} $\Delta_p$ reads on a given grid $\Gamma$ which contains an affine basis of $\R^d$  (or, equivalently, whose convex hull has a non-empty interior):
\[
\Delta_p(\xi; \Gamma) = \inf_{\lambda} \biggl\{ \Bigl( \sum_{i = 1}^N \lambda_i
	\enorm{\xi - x_i}^p \Bigr)^{1/p}: (\lambda_i)_{i=1:N} \in [0,1] ^{N}\text{ and } \sum_{i =
	1}^N \lambda_i x_i = \xi, \sum_{i = 1}^N\!\! \lambda_i = 1  \biggr\}.
\]

When $p=2$ (quadratic case),  one has the following result about the zones where the infimum $\Delta_p(\xi; \Gamma) $ is attained:   if the grid $\Gamma\subset \R^d$   contains an affine basis with its points are  in {\it general position}~--~none of its  subset of size $d+1$ lies on the same sphere~--~then it admits a unique Delaunay triangulation in the following sense (see~\cite{Rajan} or, for our setting,   Proposition~6 and Theorem~4 in~\cite{PaWi1}):

\smallskip 
\begin{enumerate}
\item  For every $\xi \!\in {\rm conv}(\Gamma)$, there exist  a unique $I= I(\xi)\subset \{1,\ldots,N\}$ of cardinality $d+1$ such that  

\smallskip 
\begin{enumerate}
\item  $(x_i)_{i\in I}$ is an affine basis,

\smallskip 
\item ${\rm conv}\{x_i,\, i\!\in I\}\cap \{x_j, \, j\!\in I^c\} =\varnothing$ (so-called {\em Delaunay property}), 

\smallskip 
\item $\Delta_p(\xi; \Gamma)$ is attained as a minimum at an $N$-tuple $\lambda_1,\ldots, \lambda_{_N}$ satisfying the constraints with   $\lambda_i=0$ if $i\notin I$.
\end{enumerate}
\smallskip 
\item If $I= I(\xi)$ as above for some $\xi\!\in  {\rm conv}(\Gamma)$, then for every $\xi' \!\in {\rm conv}\big(x_i, \, i\!\in I(\xi)\big)$, $I(\xi')= I(\xi)$.
\end{enumerate}

\smallskip A collection of simplexes $(x_i)_{i\in I}$ where $I$ is admissible for some $\xi\!\in {\rm conv}(\Gamma)$ is called a triangulation of $\Gamma$. When the points of $\Gamma$ are not in general position, several subsets $I$ of $\{1,\ldots,N\}$ can satisfy condition 1. However, if such is the case, $I$ remains admissible for all points $\xi$ in ${\rm conv}(x_i, \, i\in I)$.  Thus,  several triangulations may exist,  each one giving raise to its own splitting operator (see~\eqref{eq:ProjDel} below). A typical example is a rectangle split by one of its two diagonals which yields two triangulations, one for each diagonal.

\begin{figure}[h!]
 \centering  
  \begin{minipage}[h]{0.475\textwidth}
   \centering
   \includegraphics[height= 5cm, width= 8cm]{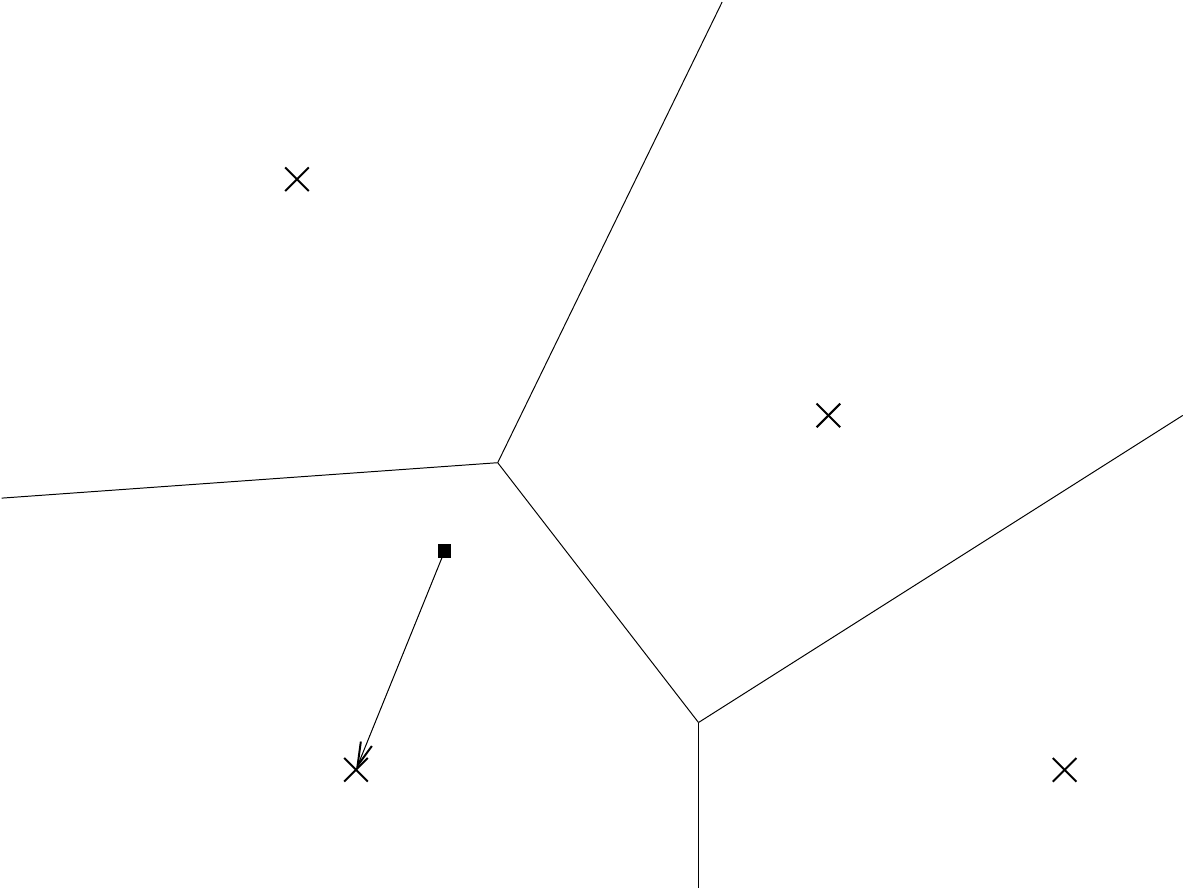} 
   \end{minipage}
   \hfill
   \begin{minipage}[h]{0.475\textwidth}
   \centering
    \includegraphics[height = 5cm, width= 8cm]{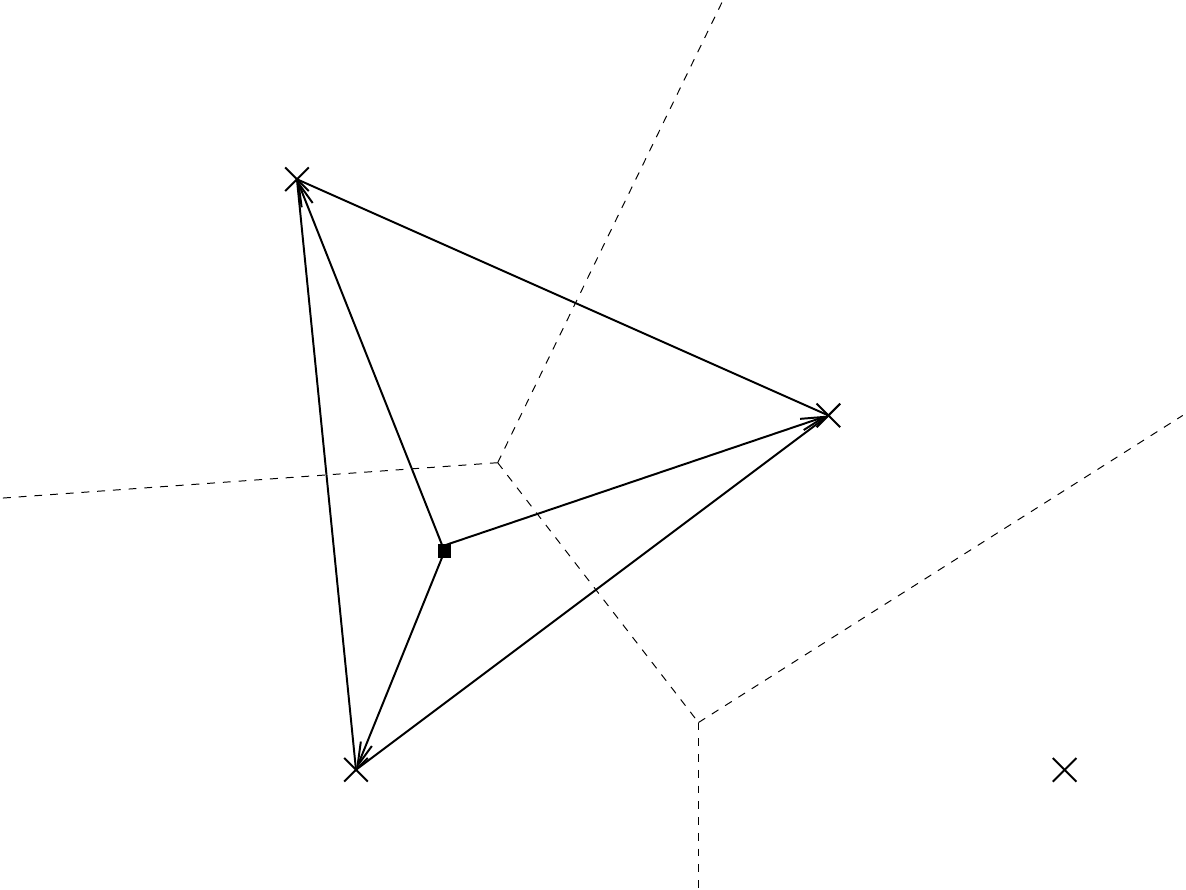}
   \end{minipage}
      \caption{\it Voronoi (left) and Delaunay (right) projections for the realization  $X(\omega)= \, ^{_\blacksquare}$.}
     \label{fig:mappings}
  \end{figure}

\medskip
It was proved in~\cite{PaWi1} that for such grids, we can
construct a {\em dual quantization projection} (or {\em splitting operator}) which is the counterpart of the nearest neighbour projection for Voronoi quantization.
This operator maps the random variable $X$ randomly to the vertices of the
Delaunay $d$-simplex (``hyper-triangle")  in which $X$ falls (see Figure~\ref{fig:mappings} further on), where the probability of mapping/projecting $X$ to a given vertex $t_i$ is determined by the $i$-th barycentric coordinate of $X$ in the (non-degenerated) $d$-simplex ${\rm conv}\{t_j: j = 1, \ldots, d+1\}$. When $p\neq2$,  the notion of Delaunay ``triangulation" can still be defined although  slightly more involved (similarly, the Voronoi cells are no longer convex when $p\neq 2$). We refer again to~\cite{PaWi1} for details. 

Mathematically speaking, let $(D_k(\Gamma))_{1\leq k\leq m}$ be a Delaunay
partition of the convex hull ${\rm conv}(\Gamma)$ of $\Gamma$. Let us denote by $\lambda^k(\xi)$ the barycentric coordinates of $\xi$
in the triangle $D_k(\Gamma)$, with the convention $\lambda^k_i(\xi)=0$ if $x_i\notin D_k(\Gamma)$. We define the  {\em dual (or Delaunay) projection operator }~--~also called {\em spliting operator}~--~by 
\begin{equation}\label{eq:ProjDel}
{\rm Proj}^{del}_{\Gamma}(\xi, u) =  \sum_{k = 1}^m \Biggl[ \sum_{i=1}^{N} x_i \cdot
	\mbox{\bf 1}_{\bigl\{\sum\limits_{j=1}^{i-1} \lambda^k_j(\xi)\,\leq\, u
 <	\sum\limits_{j=1}^{i} \lambda^k_j(\xi) \bigr\}} \Biggr]
	\mbox{\bf 1}_{D_k(\Gamma)}(\xi).
\end{equation}
Note that in~\cite{PaWi1} this projection is denoted $\mathcal{J}_\Gamma^u$ (this change is motivated by notational  consistency).  
It is clear that, by construction,  
\begin{equation}
   \forall\, \xi\!\in {\rm conv}(\Gamma), \quad \int_0^1 {\rm Proj}^{del}_{\Gamma} (\xi, u) du = \xi\label{eq:statprop}.
\end{equation}
Moreover, it follows from~\eqref{eq:ProjDel}, that
\[
\Delta_p(\xi; \Gamma) = \Bigl( \E_{\lambda_{|[0,1]}} \enorm{\xi -{\rm Proj}^{del}_{\Gamma}(\xi, U)}^p\Bigr)^{1/p},
\]
where $U$ is defined  on $([0,1], {\cal B}([0,1]),
\lambda_{|[0,1]})$ with a $\mathcal{U}\bigl([0,1]\bigr)$-distributed (so that the operator ${\rm Proj}_{\Gamma}(\xi, u) $ is
defined on this exogenous space).  Then we define  (on the product probability space $(\widetilde \Omega, \widetilde {\cal A}, \widetilde \P)$) the {\it dual} (or {\it Delaunay}) {\it quantization} 
\[
	\widehat X^{\Gamma,\text{dual}}:= {\rm Proj}^{del}_\Gamma(X,U)
\]
so that
\[
	\big\|\Delta_p(X; \Gamma)\big\|_p = \big\|X - \widehat X^{\Gamma,\text{dual}}\big\|_p\quad \mbox{ and }\quad \E(\widehat X^{\Gamma,\text{dual}}\,|\,X) = X.
\]

\noindent {\bf Remark.} $L^p$-dual quantization can be extended in a canonical way to  $L^p(\P)$-integrable random vectors by defining in a proper way the splitting operator {\em outside} the convex hull of the grid $\Gamma$, to the price of loosing the dual stationarity property.

\medskip
\noindent {\bf Optimal $L^p$-dual quantizers.} It is shown in~\cite{PaWi1} that,
for every integer $N\ge d+1$, there exists at least one {\em optimal dual quantizer} $\Gamma^{(N),del}$ at level $N\ge d+1$ which achieves the infimum $d_{p,N}(X)$ and  any such optimal dual quantizer   has cardinality $N$. Furthermore, $d_{p,N}(X)\to 0$ as $N\to + \infty$. The  convergence rate of dual quantization for bounded random vectors established in~\cite{PaWi3} turns out to be quite similar to that of the primal/Voronoi quantization as emphasized below. 

\begin{Thm}[Optimal Delaunay Quantization rates, see~\cite{PaWi1}]\label{thm:zadoretpierce2} Let $X\!\in L_{\R^d}^{\infty}(\Omega,{\cal A}, \P)$ be a bounded random vector with distribution $\P_{_X}= \varphi.\lambda_d+\nu_{_X}$ where $\nu_{_X}$ is singular. Then both claims $(a)$ and $(b)$ of Zador's Theorem remain formally true  provided the constants $\widetilde J^{vor}_{d,p}$ and $ \widetilde C^{vor}_{d,\eta,p}$ are replaced for every $p\!\in (0, +\infty)$ by  
$$\widetilde J^{del}_{d,p}:=\inf_{N\ge 1} N^{\frac 1d} d_{p,N}\big(\mathcal{U}([0,1]^d)\big)\ge \widetilde J^{vor}_{d,p},  \quad\widetilde J^{del}_{1,p}=  \Big(\frac{2}{(p+1)(p+2)}\Big)^{1/p}  \quad\mbox{ and }\quad \widetilde C^{del}_{d,\eta,p}\ge \widetilde C^{vor}_{d,\eta,p}.
$$ 
%
\end{Thm}

\noindent {\bf Remarks.} $\bullet$ Note that claim~$(b)$ remains true if the support of $\P_{_X}$  spans  an affine  subspace of $\R^d$ $A_{\mu}$ with dimension $d'<d$. However, if  such is the case the rate $N^{-1/d}$ is suboptimal, the exact one one being $N^{-1/d'}$.

\smallskip
\noindent $\bullet$ Similarly to Voronoi quantization, in a one dimensional setting, and if the distribution of $X$  is strongly unimodal, then there exist a {\em unique} $L^p$-optimal dual quantizer at any level $N\ge1$ which can be obtained, at least when $p=2$,  as the limit of the counterpart of the Lloyd  algorithm (see~\cite{JP20b}).

\medskip
\noindent {\bf Voronoi {\em versus} Delaunay quantization.} To illustrate the difference between Voronoi and Delaunay quantization (in the case $d=p=2$),
we compare in Figure~\ref{fig:mappings}  the nearest neighbor projection and the dual quantization operator.

For a given grid $\Gamma\subset\R^d$, the nearest neighbor projection  ${\rm Proj}^{vor}_\Gamma$ maps $X(\omega)$ entirely onto the generator $x_i$ of the Voronoi cell $C_i(\Gamma)$ in which $X(\omega)$ falls whereas 
the Delaunay random splitting operator ${\rm Proj}^{del}_{\Gamma}$ splits up the ``weight" $1$ of $X(\omega)$ across the vertices of the Delaunay triangle in which $X(\omega)$ falls.
Since each vertex receives here a proportion according to the barycentric coordinate of the point $X(\omega)$ in that specific Delaunay triangle,
this splitting operator fulfills a backward interpolation property,
$i.e.$ $X(\omega)$ is given by a convex combination of the vertices of the Delaunay triangle.
Finally, this property also implies the intrinsic dual stationarity condition 
$$
 \E(\widehat X^{\Gamma,\text{dual}}|X) = X.
 $$ 
 Note that, by contrast with regular Voronoi quantization where~\eqref{eq:StatioVoro} holds for optimal quadratic grids,  this dual stationarity equation is satisfies by any dual quantization grid.
 
 \bigskip
\noindent{\bf Remark.} For a comparison in one dimension, we give the example of  optimal quantizations for $\mathcal{U}([0,1])$. 
Following~\cite{PaWi1}, Section 5.1, the $L^r$-optimal dual quantizer of  $\mathcal{U}([0,1])$ (does not depend on $r$ and) is given at level $N\ge 2$ by 
$\Gamma^{(N),del}  = \Big\{ \tfrac{i-1}{N-1}: i = 1, \ldots, N  \Big\}$.
On the other hand, for Voronoi quantization,  the $L^r$-optimal quantizer of  $\mathcal{U}([0,1])$ does not depend on $r$ either, is given (see~\cite{PagSpring2018}) at level $N\ge 1$ by $\Gamma^{(N),vor}  = \Big\{ \tfrac{2i-1}{2N}: i = 1, \ldots, N  \Big\}$ and therefore made up by the midpoints of the optimal Delaunay quantizer of size $N+1$.
 Such a property does not hold for general distributions  in arbitrary dimensions. 
\end{document}